\documentclass[final,leqno,onefignum,onetabnum]{siamltex1213}

\usepackage[amssymb,thinqspace,thinspace]{SIunits}
\usepackage{amsmath}
\usepackage{amssymb}
\usepackage{graphicx}
\usepackage{subfigure}

\usepackage[mathscr]{eucal}

\newcommand{\teo}[1]{\ensuremath{\underline{\mathbf{#1}}}} 
\newcommand{\mao}[1]{\ensuremath{\mathbf{#1}}} 

\newcommand{\tet}[1]{\ensuremath{\underline{\mathbf{#1}}}} 

\newcommand{\mr}[1]{\ensuremath{\mathrm{#1}}} %

\newcommand{\mbs}[1]{\ensuremath{\boldsymbol{#1}}}
\newcommand{\names}[1]{#1}

\newcommand{\grid}[1]{#1^{\grids}}

\newcommand{\gridsub}[2]{#1^{\grids}_{#2}}
\newcommand{\gridsupsub}[3]{#1^{\grids,#2}_{#3}}
\newcommand{\flu}[1]{#1^{\flus}}
\newcommand{\flusup}[2]{#1^{\flus,#2}}
\newcommand{\flusub}[2]{#1^{\flus}_{#2}}
\newcommand{\flusupsub}[3]{#1^{\flus,#2}_{#3}}
\newcommand{\str}[1]{#1^{\strs}}
\newcommand{\strsup}[2]{#1^{\strs,#2}}
\newcommand{\strsub}[2]{#1^{\strs}_{#2}}
\newcommand{\strsupsub}[3]{#1^{\strs,#2}_{#3}}

\newcommand{\grids}{\mathscr {G}}
\newcommand{\flus}{\mathscr {F}}
\newcommand{\strs}{\mathscr {S}}

\newcommand{\IntI}{\mr I}								
\newcommand{\FSII}{\FSIIs_\mr{FSI}}			
\newcommand{\FSIIs}{\Gamma}

\newcommand{\Bound}{\ensuremath{{\Gamma}}}

\newcommand{\inv}[1]{#1^{-1}}
\newcommand{\trans}[1]{#1^{\mr{T}}}
\newcommand{\itra}[1]{#1^{\mr{-T}}}

\newcommand{\pDer}[2]{\frac{\partial #1}{\partial #2}}
\newcommand{\spDer}[2]{\frac{\partial^2 #1}{\partial #2^2}}

\newcommand{\pDerText}[2]{\partial #1 / \partial #2}
\newcommand{\tDerText}[2]{\mathrm d #1 / \mathrm d #2}

\newcommand{\SobSp}[1]{\mathcal{H}^1\left(#1\right)}

\newcommand{\LSp}[1]{\mathcal L^{2}\left(#1\right)}

\newcommand{\ltwo}{\mathrm L_2}

\newcommand{\linf}{\mathrm L_\infty}

\newcommand{\TrSp}[1]{{\mathcal{S}}_{#1}}

\newcommand{\TeSp}[1]{{\mathcal{T}}_{#1}}

\newcommand{\rhsm}{\mao r}

\newcommand{\CMatS}{\mbs{\mathscr C}_{S\!F}}  
\newcommand{\CMatF}{\mbs{\mathscr C}_{FS}}    
\newcommand{\MoM}{\mbs{\mathscr M}}   
\newcommand{\MoD}{\mbs{\mathscr D}}   
\newcommand{\MoP}{\mbs{\mathscr P}}   
\newcommand{\IMat}{\mbs{\mathscr I}}  
\newcommand{\MoMaster}{{\rm{ma}}}     
\newcommand{\MoSlave}{{\rm{sl}}}      

\newcommand{\FlMat}{\mbs{\mathscr F}}
\newcommand{\StMat}{\mbs{\mathscr S}}
\newcommand{\AleMat}{\mbs{\mathscr A}}
\newcommand{\FSIMat}{\mbs{\mathscr J}}

\newcommand{\ie}{i.e.\@ }

\newcommand{\eg}{e.g.\@ } \newcommand{\cf}{cf.\@ }

\newcommand{\bsigma}{\mbs{\sigma}}

\newcommand{\beps}{\mbs{\varepsilon}}

\newcommand{\blamb}{\mbs\lambda}

\newcommand{\bvphi}{\mbs\varphi}

\newcommand{\Dt}{\Delta t \,}
\newcommand{\tifs}{a}
\newcommand{\tiff}{b}
\newcommand{\afs}{\strsub{\alpha}{f}}
\newcommand{\aff}{\flusub{\alpha}{f}}
\newcommand{\tmid}{m}

\newcommand{\FlAle}{\ensuremath{\tau}} 
\newcommand{\FlAleI}{\ensuremath{\frac{1}{\FlAle}}}

\newcommand{\viscdyn}{\mu}  

\newcommand{\ndim}{n^\mr{dim}}

\newcommand{\Lagr}{\names{Lagrange}\@ }
\newcommand{\NS}{\names{Navier-Stokes}\@ }

\newtheorem{remark}[theorem]{Remark}

\title{A temporal consistent monolithic approach to fluid-structure interaction enabling single field predictors} 

\author{Matthias Mayr\footnotemark[2] 
\and Thomas Kl\"oppel\footnotemark[3]
\and Wolfgang A. Wall\footnotemark[5]
\and Michael W. Gee\footnotemark[2]\,\,\footnotemark[6]}

\begin{document}
\maketitle

\renewcommand{\thefootnote}{\fnsymbol{footnote}}
\footnotetext[2]{Mechanics \& High Performance Computing Group, Technische Universit\"at M\"unchen, Parkring~35, D-85748 Garching bei M\"unchen, Germany}
\footnotetext[3]{Dynamore GmbH, Industriestra{\ss}e 2, D-70565 Stuttgart, Germany. This work was performed while at Institute for Computational Mechanics, Technische Universit\"at M\"unchen.}
\footnotetext[5]{Institute for Computational Mechanics, Technische Universit\"at M\"unchen, Boltzmannstra{\ss}e~15, D-85748 Garching bei M\"unchen, Germany}
\footnotetext[6]{correspondance to gee@tum.de}
\renewcommand{\thefootnote}{\arabic{footnote}}

\begin{abstract}
We present a monolithic approach to large-deformation fluid-structure interaction (FSI) problems that allows for choosing fully implicit, single-step and single-stage time integration schemes in the structure and fluid field independently, and hence is tailored to the needs of the individual field. The independent choice of time integration schemes is achieved by temporal consistent interpolation of the interface traction. To reduce computational costs, we introduce the possibility of field specific predictors in both structure and fluid field. These predictors act on the single fields only. Possible violations of the interface coupling conditions during the predictor step are dealt with within the monolithic solution procedure. \\
We present full detail of such a generalized monolithic solution procedure, which is fully consistent in its non-conforming temporal and spatial discretization. The incorporated mortar approach allows for non-matching spatial discretizations of the fluid and solid domain at the FSI interface and is fully integrated in the resulting monolithic system of equations.
The method is applied to a variety of numerical examples. There\-by, temporal convergence rates, the special role of essential boundary conditions at the fluid-structure interface, and the positive effect of predictors are demonstrated and discussed. Emphasis is put on the comparison of different time integration schemes in fluid and structure field, for what the achieved freedom of choice of  time integrators is fully exploited.
\end{abstract}



\pagestyle{myheadings}
\thispagestyle{plain}
\markboth{\textsc{M. Mayr, T. Kl\"oppel, W. A. Wall, and M. W. Gee}}{\textsc{A Temporal Consistent Monolithic Approach to FSI}}

\section{Introduction}
\label{sec:Intro}

The interaction of fluid flow with deformable structures is of great interest in science and engineering. Especially in the case of incompressible fluid flow and finite deformation solid mechanics, it becomes very challenging to solve such coupled problems computationally. This type of problem occurs frequently in real-world physics, most notably in biomechanical or biomedical engineering.

One can distinguish between two classes of solution procedures. Solution schemes that necessitate a sequence of single field solutions and the exchange of coupling information between the fields are often referred to as \emph{partitioned} schemes (see \eg~\cite{Farhat2004,LeTallec2001}). Stability issues are discussed in~\cite{Causin2005,Foerster2007a,Joosten2010,Mok2001a}. Various acceleration techniques have been proposed in~\cite{Badia2008,Badia2008b,Kuettler2008,Kuettler2009,Nobile2008}.

Opposingly, \emph{monolithic} procedures solve both the fluid and the structural equations simultaneously within one global system of nonlinear equations. For some challenging numerical problems like channels with flexible walls~\cite{Heil2004}, thin-walled structures in the human respiratory or hemodynamic system~\cite{Kuettler2010} or for balloon-type problems like human red blood cells~\cite{Kloeppel2011}, monolithic schemes outperform partitioned procedures by far in terms of computational costs or are even the only feasible schemes to address such problems. In a monolithic framework, one looses some of the modularity of partitioned schemes but might gain great improvements in robustness and performance. Detailed performance analyses and comparisons have been carried out in~\cite{Heil2008,Kuettler2010}. Preconditioners based on block-triangular approximations of the Jacobian matrix have been introduced in~\cite{Heil2004} and extended in~\cite{Muddle2012}. In~\cite{Gee2011}, efficient preconditioners based on algebraic multigrid techniques have been developed successfully. Further preconditioning strategies can be found in~\cite{Barker2010a,Crosetto2011a}. Various techniques for the coupling of fluid and structure discretizations, where nodes do not spatially coincide at the interface, have been introduced in~\cite{Bazilevs2012,Farhat1998,Kloeppel2011,Ross2008,Ross2009}.

In opposite to space-time finite element methods, e.g. \cite{Tezduyar2006}, we discretize in time with finite difference based time integration schemes, namely with fully implicit, single-step, and single-stage time integration schemes such as the generalized-$\alpha$ method~\cite{Chung1993} in the structure field and the generalized-$\alpha$~\cite{Jansen2000} or the one-step-$\theta$ scheme in the fluid field.

In this contribution, a temporally consistent monolithic solution procedure for the interaction of incompressible fluid flow with deformable structures undergoing large deformations is presented. We allow freedom of choice for single-step, single-stage and fully implicit time integration schemes in the structure and fluid field such that the respective schemes can be tailored to the needs of the individual fields. Furthermore, individual single-field predictors are incorporated into the monolithic system leading to savings in computational costs. The predictor framework within the monolithic solver naturally enables inhomogeneous Dirichlet boundary conditions at the fluid-structure interface. For spatial discretization, finite elements for both fields are employed, whereby the interface discretizations do not need to be conforming. Exemplarily, we realize the coupling with a mortar approach that allows for complete condensation of Lagrange multipliers.

A key aspect of this paper is the temporal interpolation of tractions at the fluid-structure interface in order to consistently allow for free choice of the single field time integration schemes with non-matching time instances for the evaluation of the individual field's momentum equation. 

The derivations start with the governing equations of all fields and the coupling conditions. After discretization, a monolithic FSI residual is formulated and linearized in order to demonstrate how to build and implement a monolithic solver based on available single field codes. We discuss the discrete coupling conditions in detail with respect to meshtying of the non-matching grids at the interface and with respect to temporal consistent momentum equations of the fluid and structure field.

Selected numerical examples are used to demonstrate and discuss important properties of the proposed solution scheme. Through comparison to a FSI test case with analytical solution we report optimal temporal convergence rates. Further examples address the computational cost savings through predictors and show a thorough comparison of various combinations of different time integration schemes in fluid and structure field and their effect on the overall solution of the FSI problem.

This contribution is organized as follows. In~\S\ref{sec:Problem}, the mechanical problem at hand is introduced by means of the governing equations of both fluid and structure field as well as the coupling conditions at the fluid-structure interface. Spatial and temporal discretization is performed in~\S\ref{sec:disc} where special focus is put on the discrete interface coupling conditions. In~\S\ref{sec:monolithicsystem}, the global monolithic system of equations is assembled from all contributions derived in~\S\ref{sec:disc}. Furthermore, the condensation of the discrete \Lagr multipliers as well as the slave interface degrees of freedom is presented for both possible choices of slave and master side. For these two algorithmic variants, we report the final set of equations to be implemented in order to obtain a monolithic FSI solver. In~\S\ref{sec:examples}, we apply the proposed solution scheme to a variety of numerical examples. Finally, we close with some concluding remarks.
\section{Problem definition}
\label{sec:Problem}

\begin{figure}
\begin{center}
  \includegraphics[width=0.8\textwidth]{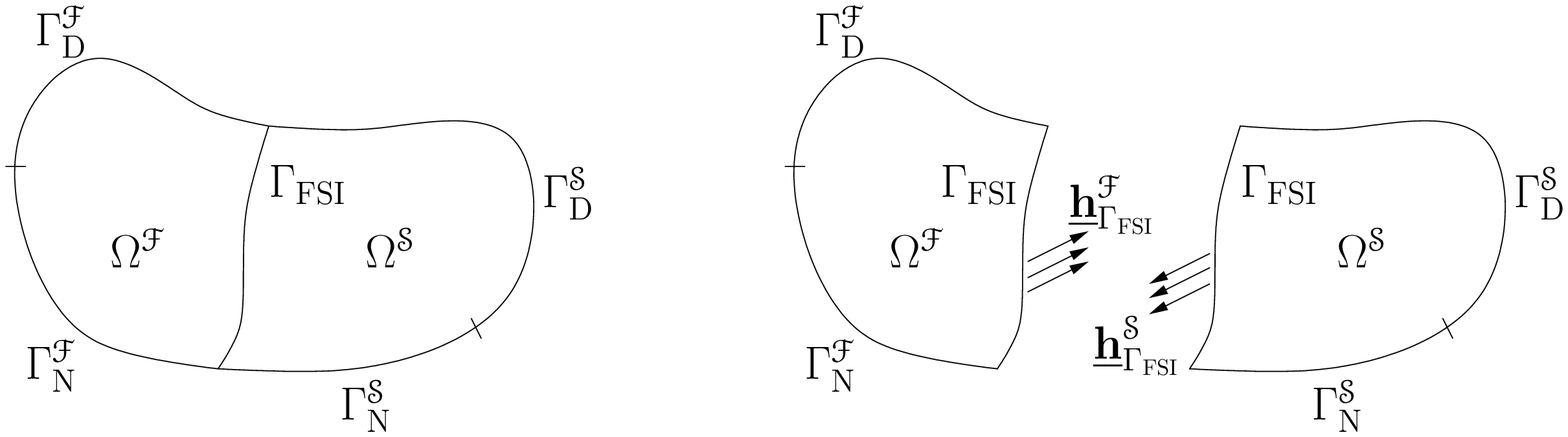}
\caption{Problem statement --- The domain~$\Omega$ is subdivided into a fluid domain~$\flu{\Omega}$ and a structural domain~$\str{\Omega}$ by the fluid-structure interface~$\FSII$. Both subdomains are bounded by Dirichlet boundaries~$\flusub{\Bound}{\mr{D}}$ and~$\strsub{\Bound}{\mr{D}}$, Neumann boundaries~$\flusub{\Bound}{\mr{N}}$ and~$\strsub{\Bound}{\mr{N}}$, and the common fluid-structure interface~$\FSII$. At the interface~$\FSII$, kinematic continuity as well as equilibrium of the interface traction fields~$\flusub{\teo h}{\FSII}$ and~$\strsub{\teo h}{\FSII}$ are required.}
\label{fig:problem}
\end{center}
\vspace{-5mm}
\end{figure}

In this section, we briefly present the governing equations for the fluid field, that is described on a deformable domain~$\flu{\Omega}$ by an \emph{Arbitrary Lagrangian Eulerian (ALE)} observer, and the structure field~$\str{\Omega}$. Both fields interact with each other at the fluid-structure interface~$\FSII$ as depicted in figure~\ref{fig:problem}, where kinematic and dynamic coupling conditions have to be satisfied. In the following, fluid quantities are denoted by the superscript~$\flu{\left(\bullet\right)}$, quantities of the ALE mesh by~$\grid{\left(\bullet\right)}$, and, finally, quantities that belong to the structure field by~$\str{\left(\bullet\right)}$. The subscript~$\left(\bullet\right)_{\FSII}$ indicates that a quantity is located at the fluid-structure interface~$\FSII$. In contrast, quantities that are located in the interior of individual field domains are marked by the subscript~$\left(\bullet\right)_{\IntI}$. To simplify notation, we often omit to state the time dependence of quantities in the sequel.

\subsection{Fluid field}
\label{ssec:FCont}

The fluid field is assumed to be governed by the instationary, incompressible \NS equations for a Newtonian fluid on a deformable domain~$\flu{\Omega}$ using an ALE description. Using an underline to indicate continuum vector or tensor valued quantities, the unknown deformation~$\grid{\teo d}(\teo x,t)$ of the deformable fluid domain~$\flu{\Omega}$ is defined by the unique mapping~$\teo\bvphi$ given as
\begin{alignat}{2}
	\label{eq:ALEdisp}
	\grid{\teo d}(\teo x,t)
	=	\teo\bvphi\left(\gridsub{\teo d}{\Bound},\teo x,t\right) 
	& \quad\mbox{ in } & \flu{\Omega}\times(0,T)\,.
\end{alignat}
The mesh deformation in the interior of the fluid domain is calculated by a mesh moving algorithm based on the boundary deformation~$\gridsub{\teo d}{\Bound}$. Then, the domain velocity~$\grid{\teo u}\left(\teo x,t\right)$ is given by
\begin{alignat}{2}
	\label{eq:ALEvel}
 	\grid{\teo u}\left(\teo x,t\right) 
  = \pDer{\teo\bvphi\left(\gridsub{\teo d}{\Bound},\teo x,t\right)}{t} 
  & \quad\mbox{ in } & \flu{\Omega}\times(0,T)\,.
\end{alignat}
In order to prevent fluid flow across the interface, it has to match the fluid velocity~$\flusub{\teo u}{\FSII}$ at the fluid-structure interface~$\FSII$
\begin{alignat}{2}
	\label{eq:ALEBdr}
 	\flusub{\teo u}{\FSII} = \gridsub{\teo u}{\FSII} = \pDer{\gridsub{\teo d}{\FSII}}{t}
  & \quad\mbox{ on } & \FSII\times(0,T)\,.
\end{alignat}
The velocity of the fluid relative to the moving background mesh is given by the ALE convective velocity~$\teo c = \flu{\teo u} - \grid{\teo u}$. Using the ALE time derivative, the incompressible \NS equations governing the fluid field on a deforming domain then read
\begin{subequations}
\label{eq:NSeq}
\begin{align}
  \label{eq:NSAle}
  \flu{\rho}\pDer{\flu{\teo u}}{t}
	+	\flu{\rho}\teo c\cdot\boldsymbol{\nabla}\flu{\teo u} 
	-	2\flu{\viscdyn}\boldsymbol{\nabla}\cdot\teo{\beps}\left(\flu{\teo u}\right)
	+	\nabla \flu p 
  & =	\flu{\rho}\flu{\teo b},\\
  \label{eq:ALEdivu}
  \boldsymbol{\nabla}\cdot\flu{\teo u} 
  & =	0\,,
\end{align}
\end{subequations}
both valid in~$\flu{\Omega}\times(0,T)$, where fluid velocity~$\flu{\teo u}$ and dynamic fluid pressure~$\flu p$ are unknown. The body force is denoted by~$\flu{\teo b}$, the strain rate tensor by~$\teo{\beps}\left(\flu{\teo u}\right) = \frac{1}{2}\left(\boldsymbol{\nabla}\flu{\teo u} + \trans{\left(\boldsymbol{\nabla}\flu{\teo u}\right)}\right)$, and the constant dynamic viscosity by~$\flu{\viscdyn}$, respectively. The fluid density~$\flu{\rho}$ is assumed to be constant. 

Given velocities~$\bar{\teo u}$ are prescribed at the Dirichlet boundary~$\flusub{\Bound}{\mr D}$. At the Neumann boundary~$\flusub{\Bound}{\mr N}$ the fluid domain is loaded with external tractions~$\flu{\bar{\teo h}}$. Additional tractions~$\flusub{\teo h}{\FSII}$ arising from the fluid-structure coupling act onto the interface portion~$\FSII$ of the boundary of the fluid subdomain~$\flu{\Omega}$. These boundary conditions read
\begin{subequations}
\label{eq:BCsFluid}
\begin{alignat}{2}
	\label{eq:BCsFluidD}
  \flu{\teo u} 
  & = \bar{\teo u} & \quad\mbox{ on }&\flusub{\Bound}{\mr D}\times(0,T)\,,\\
  \label{eq:BCsFluidN}
	\flu{\teo{\bsigma}}\cdot \flu{\teo n} 
  & = \flu{\bar{\teo h}} & \quad\mbox{ on }&\flusub{\Bound}{\mr N}\times(0,T)\,,\\
  \label{eq:BCsFluidInt}
  \flu{\teo{\bsigma}}\cdot \flu{\teo n} 
  & = \flusub{\teo h}{\FSII} & \quad\mbox{ on }&\FSII\times(0,T)\,,
\end{alignat}
\end{subequations}
where the Cauchy stress tensor~$\flu{\teo{\bsigma}}$ is defined as $\flu{\teo{\bsigma}} = -\flu{p}\tet I + 2 \flu{\viscdyn} \teo{\beps}\left(\flu{\teo u}\right)$ with the second order identity tensor~$\tet I$. The role of the interface traction~$\flusub{\teo h}{\FSII}$ in~\eqref{eq:BCsFluidInt} will be detailed in~\S\ref{ssec:FSICont}, when the coupling conditions will be discussed. As initial condition, a divergence free velocity field~$\flu{\teo u}\left(\teo x,0\right) = \flusub{\teo u}{0}\left(\teo x\right)$ with $\boldsymbol{\nabla}\cdot\flusub{\teo u}{0}\left(\teo x\right) = 0$ for $\teo x\in\flu{\Omega}$ has to be given.

Testing these equations with test functions~$\delta\flu{\teo u}$ for the momentum equation~\eqref{eq:NSAle} and~$\delta\flu p$ for the continuity equation~\eqref{eq:ALEdivu} and subsequent integration by parts gives rise to the weak form
\begin{align}
  \begin{split}
	\label{eq:ALEwf}
 	0
	= &	\left(\delta \flu{\teo u}, \flu{\rho}\pDer{\flu{\teo u}}{t}\right)_{\flu{\Omega}}
      + \left(\delta \flu{\teo u}, \flu{\rho}\teo c\cdot\mbs{\nabla}\flu{\teo u}\right)_{\flu{\Omega}}
      - \left(\mbs{\nabla}\cdot\delta \flu{\teo u}, \flu{p}\right)_{\flu{\Omega}}\\
    & + \left(\mbs{\nabla}\delta \flu{\teo u}, 2\flu{\viscdyn}\teo{\beps}(\flu{\teo u})\right)_{\flu{\Omega}}  
      - \left(\delta \flu{p},\mbs{\nabla}\cdot\flu{\teo u}\right)_{\flu{\Omega}}
      -	\left(\delta \flu{\teo u}, \flu{\rho}\flu{\teo b}\right)_{\flu{\Omega}}\\	
	& -	\left(\delta \flu{\teo u}, \flu{\bar{\teo h}}\right)_{\flusub{\Bound}{\mr N}}
	  -	\delta \flusub{W}{\FSII}\,,
  \end{split}
\end{align}
where the term
\begin{align}
  \delta \flusub{W}{\FSII}
  & = \left(\delta \teo{u},\flusub{\teo h}{\FSII}\right)_{\FSII}
\end{align}
accounts for the interface coupling and will be discussed in detail in~\S\ref{ssec:FSICont}.

\subsection{Structure field}
\label{ssec:SCont}

Without loss of generality, the structure is assumed to have a nonlinear elastic behavior. The dynamic equilibrium of forces of inertia, internal forces, and an external body force~$\str{\teo b}$ per unit undeformed volume in the undeformed structural domain~$\str{\Omega}$ is given by the nonlinear elastodynamics equation
\begin{alignat}{2}
	\label{eq:StruEqui}
	\str{\rho}\frac{\mr{d}^2 \str{\teo d}}{\mr{d}t^2}
	&	=	\mbs{\nabla}\cdot\left(\tet F\,\tet S\right)
		+	\str{\rho}\str{\teo b}
		& \quad\mbox{ in } & \str{\Omega}\times(0,T)
\end{alignat}
with the structural displacement field~$\str{\teo d}$ as the primary unknown. The structural density is denoted as~$\str{\rho}$. The internal forces are expressed in terms of the deformation gradient~$\tet F$ and the second Piola-Kirchhoff stress tensor~$\tet S$.
For the sake of simplicity we restrict ourselves to a hyperelastic material behavior. The second Piola-Kirchhoff stress tensor~$\tet S$ is then defined as $\tet S=2\cdot\pDerText{\varPsi}{\tet C}$, using the strain energy function~$\varPsi$ and the right Cauchy-Green tensor~$\tet C = \trans{\tet F}\tet F$.

The traction~$\strsub{\teo{h}}{\FSII}$ acts onto the interface portion~$\FSII$ of the boundary of the structural subdomain~$\str{\Omega}$ as shown in figure~\ref{fig:problem}. Proper Dirichlet and Neumann boundary conditions have to be prescribed on~$\strsub{\Bound}{\mr{D}}$ and~$\strsub{\Bound}{\mr{N}}$, respectively, reading
\begin{subequations}
\label{eq:BCsStr}
\begin{alignat}{2}
	\label{eq:BCsStrD}
  \str{\teo d} 
	&	= \str{\bar{\teo d}} & \quad\mbox{ on } & \strsub{\Bound}{\mr{D}}\times(0,T)\,,\\
  \label{eq:BCsStrN}
	\left(\tet F\,\tet S\right)\cdot \teo N  
	&	= \str{\bar{\teo h}} & \quad\mbox{ on } & \strsub{\Bound}{\mr{N}}\times(0,T)\,,\\
  \label{eq:BCsStrInt}
  \left(\tet F\,\tet S\right)\cdot \teo N  
  & = \strsub{\teo{h}}{\FSII} & \quad\mbox{ on } & \FSII\times(0,T)\,.
\end{alignat}
\end{subequations}
In addition, initial conditions~$\str{\teo d}\left(\teo x, 0\right) = \strsub{\teo d}{0}\left(\teo x\right)$ and~$\tDerText{\str{\teo d}}{t}\left(\teo x, 0\right) = \strsub{\dot{\teo d}}{0}\left(\teo x\right)$ have to be satisfied for~$\teo x\in\str{\Omega}$ for given initial displacement and velocity fields~$\strsub{\teo d}{0}\left(\teo x\right)$ and $\strsub{\dot{\teo d}}{0}\left(\teo x\right)$, respectively.

By multiplication of~\eqref{eq:StruEqui} with virtual displacements~$\str{\delta\teo d}$ and subsequent integration by parts one obtains the weak form
\begin{align}
\label{eq:Struwf}
  \begin{split}
    & \left(\str{\delta \teo d}, 
            \str{\rho}\frac{\str{\mr{d}^2 \teo d}}{\mr{d}t^2}\right)_{\str{\Omega}}
      +	\left(\mbs \nabla \str{\delta \teo d},\tet F\,\tet S \right)_{\str{\Omega}}
      -	\left(\str{\delta \teo d}, \str{\rho}\str{\teo b}\right)_{\str{\Omega}}\\
    & -	\left(\str{\delta \teo d}, \str{\bar{\teo h}}\right)_{\strsub{\Bound}{N}}
      -	\strsub{\delta W}{\FSII}
      = 0
  \end{split}
\end{align}
as the starting point for the finite element discretization. The term
\begin{align}
  \strsub{\delta W}{\FSII}
  & = \left(\delta \teo{d},\strsub{\teo h}{\FSII}\right)_{\FSII}
\end{align}
accounts for the influence of the interface coupling and will be discussed in the following subsection.

\subsection{Fluid-structure interface}
\label{ssec:FSICont}

Fluid field and structure field are coupled through enforcing kinematic and dynamic continuity conditions at the fluid-structure interface~$\FSII$. Physically motivated, the no-slip condition~\eqref{eq:FSInoslip_vel} is assumed that prohibits fluid flow across the fluid-structure interface and relative tangential movement of fluid and structure at the fluid-structure interface. It couples the physical fields, \ie fluid velocity field and structural displacement field. From~\eqref{eq:ALEBdr} one knows that fluid velocity  and grid velocity coincide at the fluid-structure interface, yielding~\eqref{eq:FSInoslip_dis}. Integration with respect to time finally leads to the equivalent coupling condition~\eqref{eq:StrAleCoupling}. Finally, dynamic equilibrium of interface tractions is stated in~\eqref{eq:FSIforceeq}. The coupling conditions are expressed as
\begin{subequations}
\begin{align}
  \label{eq:FSInoslip_vel}
  \pDer{\strsub{\teo d}{\FSII}}{t} & = \flusub{\teo u}{\FSII},\\
  \label{eq:FSInoslip_dis}
  \pDer{\strsub{\teo d}{\FSII}}{t} & = \pDer{\gridsub{\teo d}{\FSII}}{t},\\
  \label{eq:StrAleCoupling}
  \strsub{\teo d}{\FSII} & = \gridsub{\teo d}{\FSII},\\
  \label{eq:FSIforceeq}
  \strsub{\teo h}{\FSII} & = - \flusub{\teo h}{\FSII},
\end{align}
\end{subequations}
all valid on~$\FSII \times(0,T)$. Traction vectors~$\strsub{\teo h}{\FSII}$ and~$\flusub{\teo h}{\FSII}$ denote the traction at the fluid-structure interface onto structure and fluid field, respectively (see figure~\ref{fig:problem}).

\begin{remark}
We note that all three conditions~\eqref{eq:FSInoslip_vel}, \eqref{eq:FSInoslip_dis}, and~\eqref{eq:StrAleCoupling} are totally equivalent in the continuous regime. After temporal discretization they might differ, depending on the choice of time integration schemes in fluid and structure field.
\end{remark}

Kinematic continuity~\eqref{eq:StrAleCoupling} is imposed by a \Lagr multiplier field~$\teo\blamb$ introducing an additional field of unknowns in the coupled FSI problem. After multiplication of~\eqref{eq:StrAleCoupling} with the corresponding test function~$\delta \teo\blamb$ and a subsequent integration over the fluid-structure interface~$\FSII$ one obtains the weak form of the kinematic constraint
\begin{align}
  \label{eq:FSIweakformkinem}
  \left(\delta \teo\blamb,\strsub{\teo d}{\FSII} - \gridsub{\teo d}{\FSII} \right)_{\FSII} 
  = 0\,.
\end{align}
By identifying the interface traction~$\strsub{\teo h}{\FSII}$ onto the structure field with the \Lagr multiplier field~$\teo\blamb$ and using the dynamic equilibrium~\eqref{eq:FSIforceeq}, we can specify the interface coupling contributions~$\delta\flusub{W}{\FSII}$ and~$\delta\strsub{W}{\FSII}$ in~\eqref{eq:ALEwf} and~\eqref{eq:Struwf} as
\begin{align}
\label{eq:FSIweakform}
  \delta\flusub{W}{\FSII}
  = - \left(\delta\flu{\teo u},\teo\blamb\right)_{\FSII},
  \quad
  \delta\strsub{W}{\FSII}
  = \left(\delta\str{\teo d}, \teo\blamb\right)_{\FSII}.
\end{align}

\subsection{Weak form of coupled FSI system}
\label{ssec:WeakFSISystem}

We define the following solution spaces:
\begin{subequations}
\begin{align}
  \TrSp{\str{\teo d}} 
  & := \left\{\str{\teo d}\in\SobSp{\str\Omega}~|~\str{\teo d}=\str{\bar{\teo d}}
  \text{ on }\strsub{\Bound}{D}\right\}\\
  \TrSp{\flu{\teo u}} 
  & := \left\{\flu{\teo u}\in\SobSp{\flu\Omega}~|~\flu{\teo u}=\flu{\bar{\teo u}}
  \text{ on }\flusub{\Bound}{D}\right\}\\
  \TrSp{\flu{p}} 
  & := \left\{\flu{p}\in\LSp{\flu\Omega}\right\}\\
  \TrSp{\grid{\teo d}} 
  & := \left\{\grid{\teo d}\in\SobSp{\grid\Omega}~|~\grid{\teo d}=\grid{\bar{\teo d}}
  \text{ on }\gridsub{\Bound}{D}\right\}\\
  \TrSp{\teo\blamb}
  & := \left\{\teo\blamb\in\mathcal{H}^{-\frac{1}{2}}\left(\FSII\right)\right\}.
\end{align}
\end{subequations}
The test function spaces~$\TeSp{\str{\teo d}}$, $\TeSp{\flu{\teo u}}$, $\TeSp{\flu p}$, $\TeSp{\grid{\teo d}}$ and~$\TeSp{\teo\blamb}$ are defined as the corresponding spaces with homogeneous Dirichlet boundaries.

We finally state the overall weak problem as a combination of the weak forms~\eqref{eq:ALEwf}, \eqref{eq:Struwf}, and~\eqref{eq:FSIweakformkinem}: Find~$\str{\teo d}\in\TrSp{\str{\teo d}}$, $\flu{\teo u}\in\TrSp{\flu{\teo u}}$, $\flu p\in\TrSp{\flu p}$, $\grid{\teo d}\in\TrSp{\grid{\teo d}}$ and~$\teo\blamb\in\TrSp{\teo\blamb}$ such that
\begin{subequations}
\label{eq:FSIWF}
\begin{align}
  \begin{split}
  0 & = 
    \left(\delta \flu{\teo u}, \flu{\rho}\pDer{\flu{\teo u}}{t}\right)_{\flu{\Omega}}
    + \left(\delta \flu{\teo u}, \flu{\rho}\teo c\cdot\mbs{\nabla}\flu{\teo u}\right)_{\flu{\Omega}}
    - \left(\mbs{\nabla}\cdot\delta \flu{\teo u}, \flu{p}\right)_{\flu{\Omega}}\\
  & + \left(\mbs{\nabla}\delta \flu{\teo u}, 2\flu{\viscdyn}\teo{\beps}(\flu{\teo u})\right)_{\flu{\Omega}}      
    - \left(\delta \flu{p},\mbs{\nabla}\cdot\flu{\teo u}\right)_{\flu{\Omega}}
    - \left(\delta \flu{\teo u}, \flu{\rho}\flu{\teo b}\right)_{\flu{\Omega}}\\
  & - \left(\delta \flu{\teo u}, \flu{\bar{\teo h}}\right)_{\flusub{\Bound}{\mr N}}
    + \left(\delta\flu{\teo u},\teo\blamb\right)_{\FSII},
  \label{eq:FSIFluidWF}
  \end{split}\\
  \begin{split}
  0 & =
    \left(\str{\delta \teo d}, 
          \str{\rho}\frac{\str{\mr{d}^2 \teo d}}{\mr{d}t^2}\right)_{\str{\Omega}}
    + \left(\mbs \nabla \str{\delta \teo d},\tet F\tet S \right)_{\str{\Omega}}
    - \left(\str{\delta \teo d}, \str{\rho}\str{\teo b}\right)_{\str{\Omega}}\\
  & - \left(\str{\delta \teo d}, \str{\bar{\teo h}}\right)_{\strsub{\Bound}{N}}
    - \left(\delta\str{\teo d}, \teo\blamb\right)_{\FSII},
  \label{eq:FSIStruWF}\\
  \end{split}\\
  0 & = 
    \left(\delta \teo\blamb,
          \strsub{\teo d}{\FSII} - \gridsub{\teo d}{\FSII} \right)_{\FSII}
  \label{eq:FSICouplWF}
\end{align}
\end{subequations}
for all~$\delta\str{\teo d}\in\TeSp{\str{\teo d}}$, $\delta\flu{\teo u}\in\TeSp{\flu{\teo u}}$, $\delta\flu p\in\TeSp{\flu p}$, $\delta\grid{\teo d}\in\TeSp{\grid{\teo d}}$ and~$\delta\teo\blamb\in\TeSp{\teo\blamb}$.

\section{Discretization and mortar coupling}
\label{sec:disc}

The weak form~\eqref{eq:FSIWF} has to be discretized in space and time. For the monolithic approach presented here, the spatial discretization for fluid, ALE, and structure field is done with finite elements. The constraints at the interface are enforced using a \emph{dual} mortar method where the nodes of the fluid and structure mesh do not have to match at the interface~\cite{Kloeppel2011}. This results in a great freedom during mesh generation to tailor the meshes to the needs of the individual fields. 

For temporal discretization, fully implicit, single-step, and single-stage time integration schemes are used for all fields. Depending on the actual choices of time integration schemes and their parameters the dynamic equilibrium is formulated at an intermediate time instant~$t^{\tmid}\in\left]t^n, t^{n+1}\right]$. In general, the actual time instants for equilibrium in the fluid and structure field do not coincide, \ie~$\flusup{t}{\tmid}\neq\strsup{t}{\tmid}$. A main contribution of this work is the freedom of choosing the time integration schemes for the fluid and structure field independently and still maintaining temporal consistency between both fields, which is shown in~\S\ref{ssec:IntDisc}.

Another new aspect in time integration is that field specific predictors are allowed within the monolithic FSI framework. Due to possible predictors in the structure and fluid field, the solution at the beginning of the nonlinear iteration loop differs from the converged solution of the previous time step by additional increments
\begin{align}
  \strsupsub{\mao d}{n+1}{0}
  = \strsup{\mao d}{n} + \Delta\strsub{\mao d}{p},
  \quad
  \flusupsub{\mao u}{n+1}{0}
  = \flusup{\mao u}{n} + \Delta\flusub{\mao u}{p}
\end{align}
with the subscript~$(\bullet)_p$ indicating the predictor step. Within the predictors, both fields can evolve independently, leading to a possible violation of the kinematic continuity requirement at the fluid-structure interface, \ie possibly incompatible initial guesses for structure and fluid field. This violation can be measured and will be accounted for, when the discrete kinematic coupling conditions are derived in~\S\ref{sssec:StrAleCoupling}. Without any predictor, these additional increments vanish, \ie~$\strsub{\Delta \mao d}{p}=\mao 0$ and~$\flusub{\Delta \mao u}{p}=\mao 0$.

Discretization of the weak form~\eqref{eq:FSIWF} of the coupled FSI problem can be performed in a separated manner. Discretization of the fluid contribution~\eqref{eq:FSIFluidWF} results in the fluid residual
\begin{align}
  \flu{\rhsm}
  & = \flusub{\rhsm}{\flu{\mao u}}
    + \flusub{\rhsm}{\mao{\blamb}}
    = \left[\begin{array}{c}
        \flusub{\rhsm}{I}\\
        \flusub{\rhsm}{\FSII}\\
        \gridsub{\rhsm}{I}\\
        \gridsub{\rhsm}{\FSII}
      \end{array}\right]
    + \left[\begin{array}{c}
        \mao 0\\
        \flusub{\rhsm}{\mao\blamb,\FSII}\\
        \mao 0\\
        \mao 0
      \end{array}\right]
\end{align}
where the first term on the right hand side contains the standard fluid residual and only the second term accounts for the coupling of the fluid and structure field. Accordingly, the discretization of the structural contribution~\eqref{eq:FSIStruWF} is written as
\begin{align}
  \str{\rhsm}
  & = \strsub{\mao r}{\str{\mao d}}
    + \strsub{\rhsm}{\mao{\blamb}}
    = \left[\begin{array}{c}
        \strsub{\rhsm}{I}\\
        \strsub{\rhsm}{\FSII}
      \end{array}\right]
    + \left[\begin{array}{c}
        \mao 0\\
        \strsub{\rhsm}{\mao\blamb,\FSII}
      \end{array}\right]
\end{align}
with the first term accounting for the pure structural problem and the second term again being responsible for the fluid-structure coupling.
Finally, the weak coupling condition~\eqref{eq:FSICouplWF} is discretized yielding a residual contribution~$\rhsm^{coupl}$. Using these single field residuals, one obtains the solution of the nonlinear coupled FSI problem by solving for
\begin{align}
\label{eq:FSIResidual}
  \rhsm^{FSI}
  & = \left[\begin{array}{c}
        \str{\rhsm}\\
        \flu{\rhsm}\\
        \rhsm^{coupl}
      \end{array}\right]
    = \mao 0\,,
\end{align}
where the residual~$\rhsm^{FSI}$ depends on the structural unknowns, the fluid unknowns, and the unknown \Lagr multipliers. To solve~\eqref{eq:FSIResidual}, a Newton-type method is applied requiring the full linearization of~$\rhsm^{FSI}$ and, thus, of all single field residuals. After summarizing all unknowns of the structure field in~$\str{\mao x}$ and those of the fluid field in~$\flu{\mao x}$, respectively, the resulting linear system in Newton iteration step~$i\geq 0$ reads
\begin{align}
\label{eq:FSILinSysSchematic}
  \left[\begin{array}{ccc}
    \pDer{\strsub{\rhsm}{\str{\mao d}}}{\str{\mao x}} & \mao 0 & \pDer{\strsub{\rhsm}{\mao\blamb}}{\mao\blamb}\\
    \mao 0 & \pDer{\flusub{\rhsm}{\flu{\mao u}}}{\flu{\mao x}} & \pDer{\flusub{\rhsm}{\mao\blamb}}{\mao\blamb}\\
    \pDer{\rhsm^{coupl}}{\strsub{\mao d}{\FSIIs}} & \pDer{\rhsm^{coupl}}{\flusub{\mao u}{\FSIIs}} &\mao 0
  \end{array}\right]_{i}^{n+1}
  \left[\begin{array}{c}
    \Delta\str{\mao x}\\
    \Delta\flu{\mao x}\\
    \Delta\mao\blamb
  \end{array}\right]_{i+1}^{n+1}
  = - \left[\begin{array}{c}
        \str{\rhsm}\\
        \flu{\rhsm}\\
        \rhsm^{coupl}
      \end{array}\right]_{i}^{n+1}\,,
\end{align}
where the subscript~$\left(\bullet\right)_{\FSII}$ denoting the fluid-structure interface has been replaced by~$\left(\bullet\right)_{\FSIIs}$ to shorten the notation.
In~\eqref{eq:FSILinSysSchematic}, the splitting into degrees of freedom that belong to the interior of~$\str{\Omega}$ or~$\flu{\Omega}$ and those located at the fluid-structure interface~$\FSII$ is omitted for clarity of presentation. It will be re-introduced when the single field contributions to~\eqref{eq:FSILinSysSchematic} will be derived in the following subsections. The matrix contribution~$\pDerText{\strsub{\rhsm}{\str{\mao d}}}{\str{\mao x}}$ will be discussed in detail in~\S\ref{ssec:SDisc}. Subsection~\S\ref{ssec:FDisc} deals with the fluid discretization and will specify the matrix contribution~$\pDerText{\flusub{\rhsm}{\flu{\mao u}}}{\flu{\mao x}}$. The remaining matrix contributions that are related to the interface coupling will be addressed in~\S\ref{ssec:IntDisc}.

After solving the linear system~\eqref{eq:FSILinSysSchematic}, the update procedure is
\begin{align}
      \left[\begin{array}{c}
        \str{\mao x}\\
        \flu{\mao x}\\
        \mao\blamb
      \end{array}\right]_{i+1}^{n+1}
  & = \left[\begin{array}{c}
        \str{\mao x}\\
        \flu{\mao x}\\
        \mao\blamb
      \end{array}\right]_{i}^{n+1}
    + \left[\begin{array}{c}
        \Delta\str{\mao x}\\
        \Delta\flu{\mao x}\\
        \Delta\mao\blamb
      \end{array}\right]_{i+1}^{n+1}.
\end{align}
We stress that due to possible predictors in the single fields $\mao x_0^{n+1} \neq \mao x^{n}$.

In order to obtain the full linearization of the coupled FSI problem, we first briefly present the time discretization and linearization of the fluid, ALE, and structure field equations. A brief introduction to the mortar method will be given in~\S\ref{ssec:LagrDisc}. Afterwards, the coupling at the interface via the dual mortar method is illustrated. Furthermore, temporal consistent coupling of fluid and structure field is introduced. The assembly of the global monolithic system will then be shown in~\S\ref{sec:monolithicsystem}.

\subsection{Fluid field}
\label{ssec:FDisc}

Without loss of generality, stabilized equal-order interpolated finite elements are used for spatial discretization of the fluid field~\cite{Gresho2000}. Spatial discretization of ALE displacement, fluid velocity and fluid pressure field read
\begin{align}
\label{eq:FDiscAnsatz}
  \grid{\teo d} \approx \sum_{k=1}^{\flu{n}} \gridsub{N}{k} \gridsub{\mao d}{k},\quad
  \flu{\teo u} \approx \sum_{k=1}^{\flu{n}} \flusub{N}{k} \flusub{\mao u}{k},\quad
  \flu{p} \approx \sum_{k=1}^{\flu{n}} \flusub{N}{k} \flusub{p}{k}
\end{align}
with~$\flu{n}$ denoting the number of fluid nodes and~$\gridsub{N}{k}$ and~$\flusub{N}{k}$ being the finite element ansatz functions. 

Temporal discretization of the fluid field is done by one-step-$\theta$ or generalized-$\alpha$ schemes, \cf~\cite{Jansen2000}. 

In order to apply a Newton-type nonlinear solver, a linearization of the fluid residual~$\flusub{\rhsm}{\flu{\mao u}}(\flu{\mao u},\flu{\mao p},\grid{\mao d})$ has to be evaluated in every nonlinear iteration step~$i$. 
%
In order to prepare the coupling at the interface the nodal fluid velocities are separated: velocities of nodes on the interface are denoted by the vector~$\flusub{\mao u}{\FSIIs}$; the remaining velocity degrees of freedom are collected in a vector~$\flusub{\hat{\mao u}}{\IntI}$. We merge the vector~$\flu{\mao p}$ of nodal pressure values into the vector of inner velocities~$\flusub{\mao u}{I}
=\trans{\left[\flusub{\hat{\mao u}}{I},~\flu{\mao p}\right]}$ to simplify the notation without loosing any insight into further derivations. The introduced split into quantities belonging either to the interior or the fluid-structure interface of the fluid domain yields the matrix representation of the fluid tangent matrix contributions~$\FlMat_{\alpha\beta} = \pDerText{\flusub{\mao r}{\alpha}}{\flusub{\mao u}{\beta}}$ and~$\gridsub{\FlMat}{\alpha\beta} = \pDerText{\flusub{\mao r}{\alpha}}{\gridsub{\mao d}{\beta}}$ 
with $\alpha,\beta \in\{\IntI,\FSIIs\}$. In order to compute the solution increment~$\Delta\flusupsub{\mao x}{n+1}{i+1}$, the linear system
\begin{align}
	\label{eq:FluidLinSys}
 	\left[\begin{array}{cccc}
 		\FlMat_{\IntI\IntI}&\FlMat_{\IntI\FSIIs}&\gridsub{\FlMat}{\IntI\IntI}&\gridsub{\FlMat}{\IntI\FSIIs}\\
 		\FlMat_{\FSIIs\IntI}&\FlMat_{\FSIIs\FSIIs}&\gridsub{\FlMat}{\FSIIs\IntI}&\gridsub{\FlMat}{\FSIIs\FSIIs}
	\end{array}\right]_{i}^{n+1}
	\left[\begin{array}{c}
		\mbs\Delta \flusub{\mao u}{\IntI}\\
		\mbs\Delta \flusub{\mao u}{\FSIIs}\\
		\mbs\Delta \gridsub{\mao d}{\IntI}\\
		\mbs\Delta \gridsub{\mao d}{\FSIIs}
	\end{array}\right]_{i+1}^{n+1}
	=	-	\left[\begin{array}{c}
				\flusub{\mao r}{\IntI}\\
				\flusub{\mao r}{\FSIIs}
 			\end{array}\right]_{i}^{n+1}
\end{align}
has to be solved in every nonlinear iteration step~$i\ge 0$. Considering the mesh motion of the ALE mesh, we assume that discretization and linearization of~\eqref{eq:ALEdisp} result in an ALE system matrix~$\AleMat$. The linearized version of~\eqref{eq:ALEdisp} reads
\begin{align}
	\label{eq:ALELinSys}
 	\left[\begin{array}{cc}
 		\AleMat_{\IntI\IntI} & \AleMat_{\IntI\FSIIs}
	\end{array}\right]_{i}^{n+1}
 	\left[\begin{array}{c}
 		\mbs\Delta\gridsub{\mao d}{\IntI}\\
 		\mbs\Delta\gridsub{\mao d}{\FSIIs}
	\end{array}\right]_{i+1}^{n+1}
	=	-\gridsub{\rhsm}{\FSIIs}\,.
\end{align}
Note that the vectors of unknowns in~\eqref{eq:FluidLinSys} and~\eqref{eq:ALELinSys} both contain the mesh displacements and, thus, both systems can be combined to
\begin{align}
    \label{eq:FluidALELinSys}
  \left[\begin{array}{cccc}
    \FlMat_{\IntI\IntI}&\FlMat_{\IntI\FSIIs}&\gridsub{\FlMat}{\IntI\IntI}&\gridsub{\FlMat}{\IntI\FSIIs}\\
    \FlMat_{\FSIIs\IntI}&\FlMat_{\FSIIs\FSIIs}&\gridsub{\FlMat}{\FSIIs\IntI}&\gridsub{\FlMat}{\FSIIs\FSIIs}\\
    \mao 0 & \mao 0 & \AleMat_{\IntI\IntI} & \AleMat_{\IntI\FSIIs}
  \end{array}\right]_{i}^{n+1}
  \left[\begin{array}{c}
    \mbs\Delta \flusub{\mao u}{\IntI}\\
    \mbs\Delta \flusub{\mao u}{\FSIIs}\\
    \mbs\Delta \gridsub{\mao d}{\IntI}\\
    \mbs\Delta \gridsub{\mao d}{\FSIIs}
  \end{array}\right]_{i+1}^{n+1}
  = - \left[\begin{array}{c}
        \flusub{\mao r}{\IntI}\\
        \flusub{\mao r}{\FSIIs}\\
        \gridsub{\rhsm}{\FSIIs}
      \end{array}\right]_{i}^{n+1}\,.
\end{align}

Let us remember that the interface deformation~$\gridsub{\mao d}{\FSIIs}$ cannot evolve freely, but has to follow the fluid field or structure field interface motion. To close the fluid linear system~\eqref{eq:FluidALELinSys} a discrete coupling condition that relates fluid interface velocities~$\flusub{\mao u}{\FSIIs}$ to ALE interface displacements~$\gridsub{\mao d}{\FSIIs}$ is necessary. It will be discussed in detail in~\S\ref{ssec:IntDisc} where the discrete coupling conditions at the fluid-structure interface will be shown.

\subsection{Structure field}
\label{ssec:SDisc}
For spatial discretization of the structure field, finite elements are used. The spatial discretization of the displacement field reads
\begin{align}
  \label{eq:SDiscAnsatz}
  \str{\teo d}\approx \sum_{k=1}^{\str{n}} \strsub{N}{k} \strsub{\mao d}{k}
\end{align}
with~$\str{n}$ denoting the number of structural nodes and~$\strsub{N}{k}$ being the finite element ansatz functions. The actual choice of shape functions, element shape, and possible element technology is not of importance for the presented method. Regarding the application of additional element technology, we refer to the numerical examples in~\S\ref{sec:examples}.

Due to its efficiency and robustness, the generalized-$\alpha$ time integration scheme~\cite{Chung1993} is applied. Additionally, it offers second-order accuracy as well as the possibility of user-controlled numerical damping. 

Linearization of the structural residual~$\strsub{\mao r}{\str{\mao d}}(\str{\mao d})$ leads to the structural stiffness matrix~$\StMat$. Similar to the fluid discretization, the structural degrees of freedom are split into inner and interface degrees of freedom, yielding a block representation of the structural stiffness matrix~$\StMat_{\alpha\beta} = \pDerText{\strsub{\mao r}{\alpha}}{\strsub{\mao d}{\beta}}$, given~$\alpha,\beta\in\{\IntI, \FSIIs\}$. Putting these blocks together, the linear system
\begin{align}
 	\label{eq:StruLinSys}
 	\left[\begin{array}{cc}
 		\StMat_{\IntI\IntI} & \StMat_{\IntI\FSIIs}\\
 		\StMat_{\FSIIs\IntI} & \StMat_{\FSIIs\FSIIs}
	\end{array}\right]_{i}^{n+1}
 	\left[\begin{array}{c}
 		\mbs\Delta\strsub{\mao d}{\IntI}\\
 		\mbs\Delta\strsub{\mao d}{\FSIIs}
	\end{array}\right]_{i+1}^{n+1}
	= - \left[\begin{array}{c}
				\strsub{\mao r}{\IntI}\\
				\strsub{\mao r}{\FSIIs}
 			\end{array}\right]_{i}^{n+1}
 \end{align}
has to be solved in every Newton iteration~$i\ge 0$ in time step~$n+1$. 

\subsection{\Lagr multiplier field}
\label{ssec:LagrDisc}

The \Lagr multiplier field that enforces the interface coupling conditions~\eqref{eq:StrAleCoupling} and~\eqref{eq:FSIforceeq} is discretized using the \emph{dual} mortar method. A very brief introduction to some basics of the dual mortar method is given. For detailed derivations and theoretical background and analysis we refer to literature, \eg~\cite{Flemisch2007,Wohlmuth2000,Wohlmuth2001} and references therein.

In a mortar setup, one distinguishes between \emph{master} and \emph{slave} side~$\Bound^{\MoMaster}$ and~$\Bound^{\MoSlave}$ of the interface. The \Lagr multiplier field is discretized on the \emph{slave} side. Numerical integration takes also place on the \emph{slave} side of the interface and results in the mortar coupling matrices~$\MoD$, belonging to the slave side, and~$\MoM$, belonging to the master side, that will be introduced later. 


In opposite to classical \Lagr multiplier choices the discretization of the \Lagr multiplier field~$\blamb$ with the \emph{dual} mortar method is based on so-called \emph{dual} shape functions~$\varPhi_j$ leading to the discretization
\begin{align}
	\label{eq:MortarLagMult}
 	\teo{\blamb} \approx \sum_{j=1}^{n^\MoSlave} \varPhi_j \mao{\blamb}_j \, ,
\end{align}
with discrete nodal \Lagr multipliers $\mao{\blamb}_j$ and~$n^\MoSlave$ slave nodes. A careful construction of the \emph{dual} shape functions~$\varPhi_j$ ensures that the biorthogonality condition
\begin{align}
  \label{eq:MortarBiortho}
  \int_{\Bound^\MoSlave} \varPhi_j N^\MoSlave_k \mathrm{d} \Bound = \delta_{jk} \int_{\Bound^\MoSlave}  N^\MoSlave_k \mathrm{d} \Bound ,
\end{align} 
with the Kronecker delta~$\delta_{jk}$ is satisfied~\cite{Wohlmuth2000}. This condition plays a major role in the evaluation of the mortar coupling matrices~$\MoD$ and~$\MoM$ since it leads to a purely diagonal form of~$\MoD$. Hence, the inversion of~$\MoD$ will be computationally cheap enabling the condensation of the \Lagr multiplier field from the global system of equations as will be shown in~\S\ref{ssec:StruSplit} and~\S\ref{ssec:FluidSplit}. For further details on the application of dual mortar methods to fluid-structure interaction problems with non-conforming meshes at the fluid-structure interface the reader is referred to~\cite{Kloeppel2011}.

\subsection{Fluid-structure interface}
\label{ssec:IntDisc}

\begin{figure}
\begin{center}
  \includegraphics[width=0.35\textwidth]{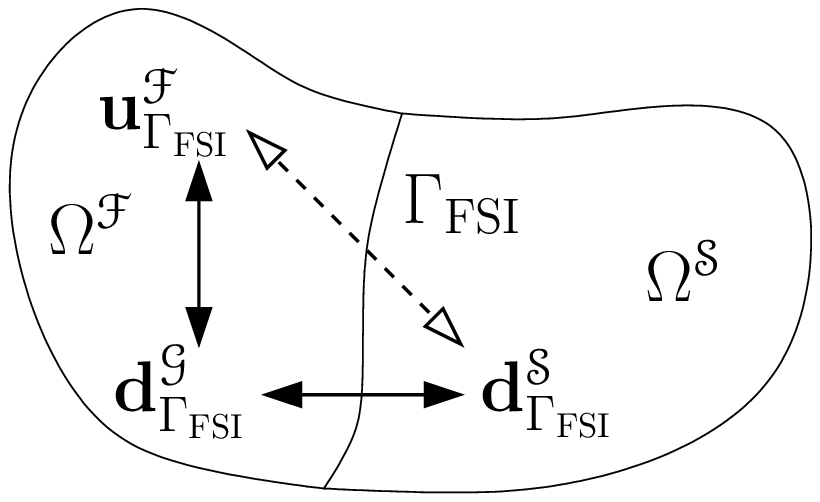}
\caption{Illustration of kinematic interface coupling conditions --- The conversion of interface fluid velocity degrees of freedom~$\flusub{\mao u}{\FSII}$ into ALE displacement degrees of freedom~$\gridsub{\mao d}{\FSII}$ happens inside the fluid field only and does not include the mortar coupling across the fluid-structure interface. The mortar coupling itself involves structure and ALE displacement degrees of freedom~$\strsub{\mao d}{\FSII}$ and~$\gridsub{\mao d}{\FSII}$, respectively. By combination of these two couplings that are illustrated by solid arrows we obtain the FSI coupling of interface fluid velocity degrees of freedom~$\flusub{\mao u}{\FSII}$ and interface structure displacement degrees of freedom~$\strsub{\mao d}{\FSII}$ that is indicated by the dashed arrow.}
\label{fig:coupling}
\end{center}
\vspace{-5mm}
\end{figure}

%
When discretizing the kinematic coupling conditions, one has to deal with the kinematic coupling of structure, fluid, and ALE degrees of freedom resulting in two separate discrete kinematic coupling conditions (see figure~\ref{fig:coupling}). On the one hand, the evolution of the fluid interface motion, described by~$\gridsub{\mao d}{\FSIIs}$, has to be related to the fluid velocity~$\flusub{\mao u}{\FSIIs}$ at the interface as already mentioned in~\S\ref{ssec:FDisc}. This coupling between fluid and ALE degrees of freedom takes place purely in the fluid domain and does not involve any structural degrees of freedom. On the other hand, a discrete version of the kinematic continuity constraint~\eqref{eq:StrAleCoupling} has to be provided for the 'meshtying' problem at the interface in order to associate  the ALE deformation~$\gridsub{\mao d}{\FSIIs}$ with the structural deformation~$\strsub{\mao d}{\FSIIs}$ at the interface. Therefore, two sets of degrees of freedom which are separated by the interface are coupled and, thus, the mortar coupling will play an important role. Both discrete kinematic coupling conditions can finally be combined to relate fluid interface velocities~$\flusub{\mao u}{\FSIIs}$ with structural interface displacements~$\strsub{\mao d}{\FSIIs}$ leading to the discrete representation of the no-slip condition~\eqref{eq:FSInoslip_vel}. The connections and dependencies of these kinematic coupling conditions are illustrated in figure~\ref{fig:coupling} and will be discussed in~\S\ref{sssec:FluAleCoupling},~\S\ref{sssec:StrAleCoupling} and~\S\ref{sssec:FluStrCoupling}.
Discretization of the dynamic constraint~\eqref{eq:FSIforceeq} will be detailed in~\S\ref{sssec:BalLinMomCoupling}.

\subsubsection{Conversion of fluid velocity and ALE displacement}
\label{sssec:FluAleCoupling}


In order to guarantee exact conservation of the volume of the fluid domain~$\flu{\Omega}$, the conversion of interface fluid velocities and interface ALE displacements has to be consistent with the fluid time integration scheme~\cite{Foerster2007}. Doing so, one can extend the geometric conservation law towards the interface. As a result of~\cite{Foerster2007}, the trapezoidal rule 
\begin{align}
  \label{eq:trapezoidal}
  \gridsupsub{\mao d}{n+1}{\FSIIs} - \gridsupsub{\mao d}{n}{\FSIIs}
  & = \frac{\Dt}{2}\left(\flusupsub{\mao u}{n+1}{\FSIIs} + \flusupsub{\mao u}{n}{\FSIIs}\right)
\end{align}
is used for the conversion of interface fluid velocities and interface ALE displacements. It is sometimes replaced by the dissipative backward Euler scheme~\cite{Foerster2007}
\begin{align}
  \label{eq:backwardeuler}
  \gridsupsub{\mao d}{n+1}{\FSIIs} - \gridsupsub{\mao d}{n}{\FSIIs}
  & = \Dt\flusupsub{\mao u}{n+1}{\FSIIs}.
\end{align}
To enable the inclusion of~\eqref{eq:trapezoidal} or~\eqref{eq:backwardeuler} into the global monolithic system, they have to be expressed in incremental form. Both can be cast into the form
\begin{align}
\label{eq:FluAleConversion}
  \Delta \gridsupsub{\mao d}{n+1}{\FSIIs,i+1}
  & = \FlAle \Delta \flusupsub{\mao u}{n+1}{\FSIIs,i+1}
    + \delta_{i0}\,\Dt \flusupsub{\mao u}{n}{\FSIIs}\,,
\end{align}
where the parameter~$\FlAle$ switches between trapezoidal rule and backward Euler scheme:
\begin{align}
  \label{eq:DefinitionTau}
  \FlAle
  & = \begin{cases}
        \frac{\Dt}{2} & \text{for trapezoidal rule~\eqref{eq:trapezoidal}}\\
        \Dt & \text{for backward Euler scheme~\eqref{eq:backwardeuler}}
      \end{cases}
\end{align}

\subsubsection{Discrete coupling condition for structural and ALE displacements}
\label{sssec:StrAleCoupling}

Inserting the spatial discretizations~\eqref{eq:FDiscAnsatz}, \eqref{eq:SDiscAnsatz}, and~\eqref{eq:MortarLagMult} for ALE displacement, structural displacement, and \Lagr multiplier field into the weak coupling condition~\eqref{eq:FSICouplWF} yields
\begin{align}
\begin{split}
  \label{eq:FSICouplDisc1}
  & \left(\sum_{j=1}^{n^\MoSlave} \varPhi_j \delta\mao{\blamb}_j,
        \sum_{k=1}^{n^{\gamma}} N^{\gamma}_{k} \strsupsub{\mao d}{\gamma}{k} - \sum_{l=1}^{n^{\varepsilon}} N^{\varepsilon}_{l} \gridsupsub{\mao d}{\varepsilon}{l} \right)_{\FSIIs} \\
  & = \sum_{j=1}^{n^\MoSlave} \delta\mao{\blamb}_j \cdot
    \left[
      \sum_{k=1}^{n^{\gamma}}\int_{\FSIIs} \varPhi_j N^{\gamma}_{k} \mr{d}\Bound \, 
      \strsupsub{\mao d}{\gamma}{k}
    - \sum_{l=1}^{n^{\varepsilon}}\int_{\FSIIs} \varPhi_j N^{\varepsilon}_{l}\mr{d}\Bound \, 
      \gridsupsub{\mao d}{\varepsilon}{l}
    \right]\\
  & = \sum_{j=1}^{n^\MoSlave} \delta\mao{\blamb}_j \cdot
      \left[\CMatS[j,k] \, \strsupsub{\mao d}{\gamma}{k} - \CMatF[j,l] \, \gridsupsub{\mao d}{\varepsilon}{l} \right]
    = 0 \quad \forall \, \delta\mao{\blamb}_j \neq \mao 0
\end{split}
\end{align}
with~$\gamma,\varepsilon\in\{\MoMaster,\MoSlave\}$, $\gamma\neq\varepsilon$ and~$n^\MoSlave$ being the number of slave nodes. Furthermore, the nodal coupling matrices
\begin{subequations}
\begin{align}
  \CMatS[j,k]
  & = \CMatS^{jk} \IMat_{\ndim}
    = \int_{\FSIIs}\varPhi_j N_{k}^{\gamma}\,\mr{d}\Bound \, \IMat_{\ndim},\,\\
  \CMatF[j,l]
  & = \CMatF^{jl} \IMat_{\ndim}
    = \int_{\FSIIs}\varPhi_j N_{l}^{\varepsilon}\,\mr{d}\Bound \, \IMat_{\ndim}\,.
\end{align}
\end{subequations}
have been introduced using an identity matrix~$\IMat_{\ndim}\in\mathbb{R}^{\ndim\times\ndim}$ with~$\ndim$ being the spatial dimension, \ie~$\ndim\in\{2,3\}$. For example, if the structure field is chosen as the slave field, \ie~$\gamma=\MoSlave$, the biorthogonality condition~\eqref{eq:MortarBiortho} can be employed to write:
\begin{align}
  \CMatS 
  & = \sum_{k=1}^{n^{\gamma}}\int_{\FSIIs} \varPhi_j N^{\MoSlave}_{k} \mr{d}\Bound
    = \sum_{k=1}^{n^{\gamma}}\delta_{jk}\int_{\FSIIs} N^{\MoSlave}_{k} \mr{d}\Bound
\end{align}
When choosing the fluid field as the slave field, the coupling matrix~$\CMatF$ takes a diagonal form in an analogous way. Full details on the numerical evaluation of the mortar integrals are given in~\cite{Puso2004,Puso2004a,Puso2004b}.

Assembling the nodal coupling matrices leads to global coupling matrices~$\CMatS$ and~$\CMatF$, which are used to formulate the kinematic coupling residual
\begin{align}
  \rhsm^{coupl}
  & = \CMatS\strsub{\mao d}{\FSIIs} - \CMatF\gridsub{\mao d}{\FSIIs}
    = \mao 0.
\end{align}
Its linearization yields the kinematic coupling constraint
\begin{align}
  \label{eq:MortarLinSys1}
  \CMatS \mbs\Delta \strsupsub{\mao d}{n+1}{\FSIIs,i+1}
  - \CMatF \mbs\Delta \gridsupsub{\mao d}{n+1}{\FSIIs,i+1}
  & = - \delta_{i0}\,\CMatS\Delta\strsub{\mao d}{\FSIIs,p} 
\end{align}
formulated in incremental form. The violation of the interface continuity requirement due to possible non-constant predictors is measured by~$\Delta\strsub{\mao d}{\FSIIs,p}$ and accounted for by the right hand side term, which is necessary only in the first nonlinear iteration step~$i=0$. Due to the linearity of the kinematic coupling condition, the kinematic interface continuity requirement is guaranteed to be satisfied for all nonlinear iteration steps~$i>0$. In the case of conforming interface discretizations, all mortar projection operators reduce to diagonal matrices with area weights on the main diagonal as well as the interface constraints collapse to the trivial case of condensable point-wise constraints. Furthermore, note that~\eqref{eq:MortarLinSys1} guarantees exact kinematic continuity in the discrete setting, even if different time integrations schemes are employed. 

\subsubsection{Discrete coupling condition for structural displacements and fluid velocities}
\label{sssec:FluStrCoupling}

With the discrete coupling conditions~\eqref{eq:FluAleConversion} and~\eqref{eq:MortarLinSys1}, all necessary conditions are at hand to assemble the global monolithic system. However, a direct conversion of fluid velocities and structural displacements at the fluid-structure interface can be derived by replacing the interface ALE displacements. On the one hand, this emphasizes the fact that the ALE field is not a physical field but rather an auxiliary field to describe the fluid motion. On the other hand, this eases the notation of the global monolithic system when it comes to choosing master and slave side in the context of the dual mortar method. 

Combining~\eqref{eq:FluAleConversion} and~\eqref{eq:MortarLinSys1} results in
\begin{align}
  \label{eq:DisVelConversion}
  \CMatS \Delta \strsupsub{\mao d}{n+1}{\FSIIs,i+1}
    + \delta_{i0}\,\CMatS \Delta \strsub{\mao d}{\FSIIs,p}
  & = \FlAle\, \CMatF \Delta \flusupsub{\mao u}{n+1}{\FSIIs,i+1}
    + \delta_{i0}\,\Dt \CMatF \flusupsub{\mao u}{n}{\FSIIs}
\end{align}
Note that~\eqref{eq:DisVelConversion} does not take the role of an additional coupling condition. It is just a redundant reformulation of~\eqref{eq:FluAleConversion} and~\eqref{eq:MortarLinSys1}.

\subsubsection{Contributions to the balances of linear momentum}
\label{sssec:BalLinMomCoupling}

Discretization of the fluid contribution~$\delta\flusub{W}{\FSII}$ in~\eqref{eq:FSIweakform} yields:
\begin{align}
  \label{eq:FSIFluCouplDisc}
  \begin{split}
    \left(
      \sum_{l=1}^{n^{\varepsilon}} N_l^\varepsilon \delta\flusupsub{\mao u}{\varepsilon}{l}, 
      \sum_{j=1}^{n^\MoSlave} \varPhi_j \mao{\blamb}_j
    \right)_{\FSIIs}
  & = \sum_{l=1}^{n^{\varepsilon}}\delta\flusupsub{\mao u}{\varepsilon}{l} \cdot
      \left[
        \sum_{j=1}^{n^\MoSlave} \int_{\FSIIs} N_l^\varepsilon \varPhi_j\,\mr{d}\Bound \,\mao{\blamb}_j
      \right]\\
  & = \sum_{l=1}^{n^{\varepsilon}}\delta\flusupsub{\mao u}{\varepsilon}{l} \cdot
      \trans{\CMatF}[j,l] \mao{\blamb}_j.
  \end{split}
\end{align}
Discretizing the structural contribution~$\delta\strsub{W}{\FSII}$ in~\eqref{eq:FSIweakform} results in
\begin{align}
  \label{eq:FSIStrCouplDisc}
  \begin{split}
  - \left(
      \sum_{k=1}^{n^{\gamma}} N_k^\gamma \delta\strsupsub{\mao d}{\gamma}{k}, 
      \sum_{j=1}^{n^\MoSlave} \varPhi_j \mao{\blamb}_j
    \right)_{\FSIIs}
  & = - \sum_{k=1}^{n^{\gamma}}\delta\strsupsub{\mao d}{\gamma}{k} \cdot
      \left[
        \sum_{j=1}^{n^\MoSlave} \int_{\FSIIs} N_k^\gamma \varPhi_j\,\mr{d}\Bound \,\mao{\blamb}_j
      \right]\\
  & = - \sum_{k=1}^{n^{\gamma}}\delta\strsupsub{\mao d}{\gamma}{k} \cdot
      \trans{\CMatS}[j,k] \mao{\blamb}_j.
  \end{split}
\end{align}
In both expressions~\eqref{eq:FSIFluCouplDisc} and~\eqref{eq:FSIStrCouplDisc} one can identify the transposes of the coupling matrices~$\CMatF$ and~$\CMatS$ which have been already introduced in~\S\ref{sssec:StrAleCoupling}. Note that the coupling matrices depend on the initial mesh configuration, only. 

In the following, we assume that time integration schemes in both fluid and structure field evaluate the single field dynamic equilibrium at intermediate time instances~$\flusup{t}{\tmid}$ and~$\strsup{t}{\tmid}$, respectively. The intermediate time instances will be indicated by the superscript~$(\bullet)^\tmid$. In general, these time instances do not coincide, i.e. $\flusup{t}{\tmid}\neq\strsup{t}{\tmid}$. Based on the discretized weak forms~\eqref{eq:FSIFluCouplDisc} and~\eqref{eq:FSIStrCouplDisc} and using linear interpolations as usual for fully implicit, single-step, single-stage time integration schemes, one can write the residual contributions~$\flusupsub{\rhsm}{\tmid}{\blamb}$ and~$\strsupsub{\rhsm}{\tmid}{\blamb}$ as
\begin{align}
\label{eq:CoupResConsistent}
  \flusupsub{\mao r}{\tmid}{\mao\blamb,i}
  = \trans{\CMatF}\left(\tiff\mao\blamb^{n}+(1-\tiff)\mao\blamb^{n+1}_i\right), 
  \quad
  \strsupsub{\mao r}{\tmid}{\mao\blamb,i}
  = - \trans{\CMatS}\left(\tifs\mao\blamb^{n}+(1-\tifs)\mao\blamb^{n+1}_i\right)
\end{align}
with time interpolation factors~$\tifs$ and~$\tiff$ chosen depending on the specific field time integrators.

\begin{remark}
  Factors~$\tifs$ and~$\tiff$ are always chosen equal to the weighting of the previous solution in~\eqref{eq:CoupResConsistent}. For example, when using generalized-$\alpha$ time integration~\cite{Chung1993} for the structure field and generalized-$\alpha$ time integration~\cite{Jansen2000} for the fluid field, time interpolation factors have to be chosen as $\tifs = \afs$ and $\tiff = 1 - \aff$.
\end{remark}

Since the weak forms~\eqref{eq:FSIweakform} are linear in the displacement field~$\str{\teo d}$ and the velocity field~$\flu{\teo u}$, the linearizations of the residual terms~\eqref{eq:CoupResConsistent} are just the coupling matrices themselves wherein the temporal interpolation factors occur, too:
\begin{align}
  \pDer{\flusupsub{\rhsm}{\tmid}{\mao\blamb,i}}{\mao{\blamb}^{n+1}_{i}} 
  = \left(1-\tiff\right) \trans{\CMatF},
  \quad
  \pDer{\strsupsub{\rhsm}{\tmid}{\mao\blamb,i}}{\mao{\blamb}^{n+1}_{i}} 
  = - \left(1-\tifs\right) \trans{\CMatS}.
\end{align}
With these linearizations, we have finally specified all contributions to the linear system~\eqref{eq:FSILinSysSchematic}.

\begin{remark}
\label{rem:InterfaceEnergy}
  We can calculate the amount of energy production per time step at the fluid-structure interface due to differences in temporal discretization of the individual fields as
  \begin{align*}
  	& \Delta E_{\FSIIs}^{n\rightarrow n+1} 
      = \strsupsub{E}{n\rightarrow n+1}{\FSIIs} + \flusupsub{E}{n\rightarrow n+1}{\FSIIs}\\
    & = \left(\tifs\blamb^n+(1-\tifs)\blamb^{n+1}\right)\left(\strsupsub{\mao{d}}{n+1}{\FSIIs}-\strsupsub{\mao{d}}{n}{\FSIIs}\right) 
      - \left(\tiff\blamb^n+(1-\tiff)\blamb^{n+1}\right)\left(\gridsupsub{\mao{d}}{n+1}{\FSIIs}-\gridsupsub{\mao{d}}{n}{\FSIIs}\right)\\
    & = \left((\tifs-\tiff)\blamb^{n}+(\tiff-\tifs)\blamb^{n+1}\right)\left(\strsupsub{\mao{d}}{n+1}{\FSIIs}-\strsupsub{\mao{d}}{n}{\FSIIs}\right)
  \end{align*}
  where the discrete kinematic coupling constraint~\eqref{eq:MortarLinSys1} has been exploited. We make the following observations:
  \begin{itemize}
    \item The energy production per step vanishes for~$\tifs-\tiff\rightarrow 0$, \ie as time instances of evaluating structure and fluid coupling tractions coincide: $\strsup{t}{\tmid}-\flusup{t}{\tmid}\rightarrow 0$.
    \item Since~$\strsupsub{\mao{d}}{n+1}{\FSIIs}-\strsupsub{\mao{d}}{n}{\FSIIs}\propto\Dt$, the energy production per step reduces as~$\Dt\rightarrow 0$. Thus, we call the scheme \emph{temporal consistent}.
  \end{itemize}
  These observations can be reproduced in numerical studies (see {\S}5.2).
\end{remark}

\section{Monolithic FSI system}
\label{sec:monolithicsystem}

We can now put all linearized single field systems together to the global monolithic linear system 
\begin{align}
  \FSIMat^{n+1}_{i}\Delta\mao x_{i+1}^{n+1} = -\rhsm_i^{FSI,n+1}
\end{align}
that has to be solved in the $i^{th}$ iteration step of the Newton-type nonlinear solution algorithm in time step~$n+1$. The Jacobian matrix of the coupled FSI system reads
\begin{subequations}
\label{eq:MFSI3System}
\begin{align}
\label{eq:MFSI3fJac}
\FSIMat^{n+1}_{i}
& = \left[\begin{array}{ccccccc}
			\StMat_{\IntI\IntI} & \StMat_{\IntI \FSIIs}\\
			\StMat_{\FSIIs\IntI} & \StMat_{\FSIIs\FSIIs} &&&&&-(1-\tifs)\trans{\CMatS}\\
			&&\FlMat_{\IntI\IntI}&\FlMat_{\IntI\FSIIs}&\gridsub{\FlMat}{\IntI\IntI}&\gridsub{\FlMat}{\IntI\FSIIs}\\
			&&\FlMat_{\FSIIs\IntI}&\FlMat_{\FSIIs\FSIIs}&\gridsub{\FlMat}{\FSIIs \IntI}&\gridsub{\FlMat}{\FSIIs\FSIIs}&(1-\tiff)\trans{\CMatF}\\
			&&&&\AleMat_{\IntI\IntI}&\AleMat_{\IntI\FSIIs}\\
			&-\CMatS&&\FlAle\CMatF
		\end{array}\right]_i^{n+1},
\end{align}
where $\mao 0$-blocks are omitted for the sake of clarity. The global solution increment vector 
\begin{align}
\label{eq:MFSI3fInc}
\trans{{\Delta\mao x_{i+1}^{n+1}}}
& =	\left[\begin{array}{ccccccc}
			\trans{{\strsub{\Delta \mao d}{\IntI}}}&
			\trans{{\strsub{\Delta \mao d}{\FSIIs}}}&
			\trans{{\flusub{\Delta \mao u}{\IntI}}}&
			\trans{{\flusub{\Delta \mao u}{\FSIIs}}}&
			\trans{{\gridsub{\Delta \mao d}{\IntI}}}&
      \trans{{\gridsub{\Delta \mao d}{\FSIIs}}}&
			\trans{{\mao\blamb}}
		\end{array}\right]_{i+1}^{n+1}
\end{align}
contains the primary unknowns of each field as well as the \Lagr multiplier field. Using the Kronecker delta~$\delta_{i0}$, the corresponding residual vector is given by
\begin{align}
\label{eq:MFSI3fRHS}
\rhsm_i^{FSI,n+1}
=	\left[\begin{array}{c}
		\strsub{\rhsm}{\IntI}\\
		\strsub{\rhsm}{\FSIIs}\\
		\flusub{\rhsm}{\IntI}\\
		\flusub{\rhsm}{\FSIIs}\\
		\gridsub{\rhsm}{\FSIIs}\\
		\mao 0
 	\end{array}\right]_i^{n+1}
+	\left[\begin{array}{c}
		\mao 0\\
		-\tifs\trans{\CMatS}\mao\blamb^n\\
		\mao 0\\
		\tiff\trans{\CMatF}\mao\blamb^n\\
		\mao 0\\
		\mao 0
 	\end{array}\right]
+	\delta_{i0}
	\left[\begin{array}{c}
		\mao 0\\
		\mao 0\\
		\mao 0\\
		\mao 0\\
		\mao 0\\
		\Dt\CMatF\flusupsub{\mao u}{n}{\FSIIs} - \CMatS\strsub{\Delta \mao d}{\FSIIs,p}.
	\end{array}\right].
\end{align}
\end{subequations}
To close system~\eqref{eq:MFSI3System} still a coupling of interface ALE degrees of freedom to the motion of the fluid-structure interface~$\FSII$ is required. The missing equation will be added after master and slave side have been chosen. Doing so, one can formulate the description of the interface motion in terms of master degrees of freedom, which will be a good starting point for condensation of \Lagr multipliers and interface slave degrees of freedom. The missing coupling equation as well as the process of condensation will be detailed in~\S\ref{ssec:StruSplit} and~\S\ref{ssec:FluidSplit}.

Due to the $\mao 0$-block on the main diagonal, the global monolithic linear system~\eqref{eq:MFSI3System} is of saddle-point type. In order to circumvent the saddle-point like system to be able to use efficient FSI specific linear solvers~\cite{Gee2011} designed for the case of conforming discretizations, the unknown \Lagr multipliers~$\mao\blamb^{n+1}$ will be condensed, yielding a problem with structural displacement, fluid velocity and pressure as well as ALE grid displacement degrees of freedom as the only unknowns. Employing the kinematic coupling conditions that were derived in~\S\ref{ssec:IntDisc} one can condense the interface degrees of freedom of the slave side from the global system of equations. Depending on the choice of master and slave side, the balance of linear momentum either of the structure or of the fluid field is used to condense the discrete \Lagr multipliers. For the process of condensation we exploit the biorthogonality property~\eqref{eq:MortarBiortho} of the dual mortar method since it allows a computationally cheap inversion of the slave side's mortar matrix~$\MoD$.

After complete condensation, the interface motion is purely described and handled in terms of unknowns of the master field. Thus, we distinguish two algorithmic variants, namely \emph{fluid-handled interface motion} and \emph{structure-handled interface motion}. In the following, the two possible choices of master and slave side are discussed and the condensation process as well as the final linear systems of equations will be shown.

\begin{remark}
  In case, one wants to use standard shape functions for the \Lagr multipliers, the condensation is numerically very costly or even unfeasable. Then, the saddle-point type system~\eqref{eq:MFSI3System} can be solved with appropriate saddle-point solvers.
\end{remark}

\subsection{Fluid-handled interface motion}
\label{ssec:StruSplit}

Let us first consider the variant where the interface motion is expressed in terms of fluid velocity degrees of freedom, \ie the fluid field is chosen as the master field and the structure field as the slave field, respectively. Since the fluid field has been chosen as master, we can identify the mortar matrices as~$\MoD=\CMatS$ and~$\MoM=\CMatF$. The coupling of the interface ALE displacement to the interface motion is expressed in terms of the master's side interface degrees of freedom, \ie in terms of interface fluid velocities. This coupling has already been stated in~\eqref{eq:FluAleConversion} and will be used to close the monolithic system of equations yielding the Jacobian matrix
\begin{subequations}
\label{eq:SystemSSfull}
\begin{align}
\label{eq:JacSSfull}
\FSIMat^{n+1}_i
& = \left[\begin{array}{ccccccc}
			\StMat_{\IntI\IntI} & \StMat_{\IntI \FSIIs}\\
			\StMat_{\FSIIs\IntI} & \StMat_{\FSIIs\FSIIs} &&&&&-(1-\tifs)\trans{\MoD}\\
			&&\FlMat_\mr{II}&\FlMat_{\IntI\FSIIs}&\gridsub{\FlMat}{\IntI\IntI}&\gridsub{\FlMat}{\IntI\FSIIs}\\
			&&\FlMat_{\FSIIs\IntI}&\FlMat_{\FSIIs\FSIIs}&\gridsub{\FlMat}{\FSIIs \IntI}&\gridsub{\FlMat}{\FSIIs\FSIIs}&(1-\tiff)\trans{\MoM}\\
			&&&&\AleMat_{\IntI\IntI}&\AleMat_{\IntI\FSIIs}\\
			&-\MoD&&\FlAle\MoM\\
      &&&\FlAle\IMat&&-\IMat
		\end{array}\right]_{i}^{n+1}
\end{align}
and the residual vector
\begin{align}
\label{eq:rhsSSfull}
\rhsm_i^{n+1}
=	\left[\begin{array}{c}
		\strsub{\rhsm}{\IntI}\\
		\strsub{\rhsm}{\FSIIs}\\
		\flusub{\rhsm}{\IntI}\\
		\flusub{\rhsm}{\FSIIs}\\
		\gridsub{\rhsm}{\FSIIs}\\
    \mao 0\\
		\mao 0
 	\end{array}\right]_i^{n+1}
+	\left[\begin{array}{c}
		\mao 0\\
		-\tifs\trans{\MoD}\mao\blamb^{n}\\
		\mao 0\\
		\tiff\trans{\MoM}\mao\blamb^{n}\\
		\mao 0\\
    \mao 0\\
		\mao 0
 	\end{array}\right]
+	\delta_{i0}
	\left[\begin{array}{c}
		\mao 0\\
		\mao 0\\
		\mao 0\\
    \mao 0\\
		\mao 0\\
		\Dt\MoM\,\flusupsub{\mao u}{n}{\FSIIs} - \MoD\strsub{\Delta \mao d}{\FSIIs,p}\\
    \Dt\flusupsub{\mao u}{n}{\FSIIs}
	\end{array}\right].
\end{align}
\end{subequations}

The sixth row that brings the discrete kinematic constraint~\eqref{eq:DisVelConversion} into the system can be resolved for the structural interface displacement increment
\begin{align}
\label{eq:KinemConstrSS}
  \strsupsub{\Delta \mao d}{n+1}{\FSIIs} 
  & = \FlAle\MoP\,\flusupsub{\Delta\mao u}{n+1}{\FSIIs}
    + \delta_{i0}\,\Dt\MoP\,\flusupsub{\mao u}{n}{\FSIIs}
    - \delta_{i0}\,\strsub{\Delta\mao d}{\FSIIs,p}
\end{align}
with the mortar projection matrix 
\begin{align}
\label{eq:MortarP}
  \MoP = \inv\MoD\MoM
\end{align}
that can be efficiently computed due to the diagonal form of~$\MoD$ (\cf~\cite{Wohlmuth2000}). The coupling of fluid velocities and ALE displacements at the fluid-structure interface is given by the last row in~\eqref{eq:SystemSSfull} and reads
\begin{align}
  \label{eq:ALEConstrSS}
  \Delta\gridsupsub{\mao d}{n+1}{\FSIIs,i+1}
  & = \FlAle\Delta\flusupsub{\mao u}{n+1}{\FSIIs,i+1}
    + \delta_{i0}\,\Dt\flusupsub{\mao u}{n}{\FSIIs}.
\end{align}
From the balance of linear momentum of the structural interface degrees of freedom together with~\eqref{eq:KinemConstrSS},  the unknown \Lagr multipliers are expressed by
\begin{align}
\label{eq:RecoverLambdaSS}
\begin{split}
	\blamb^{n+1} 
		= -\frac{\tifs}{1-\tifs}\blamb^{n}
	&	+	\frac{1}{1-\tifs}\itra{\MoD}
			\left(\strsupsub{\rhsm}{n+1}{\FSIIs}
						+	\StMat_{\FSIIs\IntI}\strsupsub{\Delta \mao d}{n+1}{\IntI,i+1}
            + \FlAle\StMat_{\FSIIs\FSIIs}\MoP\Delta\flusupsub{\mao u}{n+1}{\FSIIs,i+1}\right)\\
	&	+	\delta_{i0} \frac{1}{1-\tifs}\itra{\MoD}
      \left(\Dt\StMat_{\FSIIs\FSIIs}\MoP\flusupsub{\mao u}{n}{\FSIIs}
            - \StMat_{\FSIIs\FSIIs}\strsub{\Delta \mao d}{\FSIIs,p}\right).
\end{split}
\end{align}
Equation~\eqref{eq:RecoverLambdaSS} is used to recover the \Lagr multiplier solution at the end of each time step as a postprocessing step.

Using~\eqref{eq:KinemConstrSS}, \eqref{eq:ALEConstrSS}, and~\eqref{eq:RecoverLambdaSS}, we are able to condense the system of equations. The condensed linear system with fluid-handled interface motion consists of the Jacobian matrix
\begin{subequations}
\label{eq:SystemSScond}
\begin{align}
\label{eq:JacSScond}
\FSIMat^{n+1}_i
& = \left[\begin{array}{cccc}
			\StMat_{\IntI\IntI} && \FlAle\StMat_{\IntI\FSIIs}\MoP\\
			&\FlMat_\mr{II}&\FlMat_{\IntI\FSIIs}+\FlAle\gridsub{\FlMat}{\IntI\FSIIs}&\gridsub{\FlMat}{\mr{II}}\\
			\frac{1-\tiff}{1-\tifs}\trans{\MoP}\StMat_{\FSIIs \IntI}&\FlMat_{\FSIIs\IntI}&\FlMat_{\FSIIs\FSIIs}+\FlAle\gridsub{\FlMat}{\FSIIs\FSIIs}+\frac{1-\tiff}{1-\tifs}\FlAle\trans{\MoP}\StMat_{\FSIIs\FSIIs}\MoP&\gridsub{\FlMat}{\FSIIs\IntI}\\
			&&\FlAle\AleMat_{\IntI\FSIIs}&\AleMat_{\IntI\IntI}\\
		\end{array}\right]_{i}^{n+1},
\end{align}
the solution increment vector
\begin{align}
\label{eq:SolSScond}
\trans{{\Delta\mao x_{i+1}^{n+1}}}
& =	\left[\begin{array}{cccc}
			\trans{{\strsub{\Delta \mao d}{\IntI}}}&
			\trans{{\flusub{\Delta \mao u}{\IntI}}}&
			\trans{{\flusub{\Delta \mao u}{\FSIIs}}}&
			\trans{{\gridsub{\Delta \mao d}{\IntI}}}
		\end{array}\right]_{i+1}^{n+1},
\end{align}
and the residual vector
\begin{align}
\label{eq:rhsSScond}
  \begin{split}
  \rhsm_i^{n+1}
  &	=	\left[\begin{array}{c}
        \strsub{\rhsm}{\IntI}\\
        \flusub{\rhsm}{\IntI}\\
        \flusub{\rhsm}{\FSIIs} + \frac{1-\tiff}{1-\tifs}\trans{\MoP}\strsub{\rhsm}{\FSIIs}\\
        \gridsub{\rhsm}{\FSIIs}
      \end{array}\right]_i^{n+1}
    +	\left[\begin{array}{c}
        \mao 0\\
        \mao 0\\
        \left(\tiff-\frac{\tifs(1-\tiff)}{1-\tifs}\right)\trans{\MoM}\mao\blamb^n\\
        \mao 0
      \end{array}\right]\\
  &	+	\delta_{i0}
      \left[\begin{array}{c}
        \Dt\StMat_{\IntI \FSIIs}\MoP\flusupsub{\mao u}{n}{\FSIIs} - \StMat_{\IntI \FSIIs}\strsub{\Delta \mao d}{\FSIIs,p}\\
        \Dt\gridsub{\FlMat}{\IntI\FSIIs} \flusupsub{\mao u}{n}{\FSIIs}\\
        \Dt\gridsub{\FlMat}{\FSIIs\FSIIs} \flusupsub{\mao u}{n}{\FSIIs} + \frac{1-\tiff}{1-\tifs}\Dt\trans{\MoP}\StMat_{\FSIIs\FSIIs}\MoP\flusupsub{\mao u}{n}{\FSIIs} - \frac{1-\tiff}{1-\tifs}\trans{\MoP}\StMat_{\FSIIs\FSIIs}\strsub{\Delta \mao d}{\FSIIs,p}\\
        \Dt\AleMat_{\IntI\FSIIs} \flusupsub{\mao u}{n}{\FSIIs}
      \end{array}\right].
  \end{split}
\end{align}
\end{subequations}

\begin{remark}
\label{rem:InterfaceDBCsSS}
  When parts of the mortar interface are subject to essential boundary conditions, \cite{Puso2003} suggests to apply them only on the master side of the interface in order to avoid stability problems. For the fluid-handled interface motion, this means that at the interface only the fluid field is allowed to carry Dirichlet boundary conditions. They will be imposed on the structure field weakly via the mortar coupling. 
\end{remark}

\subsection{Structure-handled interface motion}
\label{ssec:FluidSplit}

The other possibility is to describe the interface motion in terms of structural displacements, \ie the structure field is the master field and the fluid field the slave field, respectively. It is obtained by choosing the mortar matrices as~$\MoD = \CMatF$ and~$\MoM = \CMatS$. The coupling of interface ALE degrees of freedom is still governed by~\eqref{eq:FluAleConversion}. However, the fluid interface velocities will be condensed from the global system of equations and, thus, this coupling is expressed in terms of structural displacements:
\begin{align}
  \Delta\gridsupsub{\mao d}{n+1}{\FSIIs,i+1}
  & = \MoP\Delta\strsupsub{\mao d}{n+1}{\FSIIs,i+1}
    + \delta_{i0}\,\MoP\Delta\strsub{\mao d}{\FSIIs,p}\,.
\end{align}
In this case, the global monolithic linear system consists of the Jacobian matrix
\begin{subequations}
\label{eq:SystemFSfull}
\begin{align}
\label{eq:JacFSfull}
\FSIMat^{n+1}_i
& = \left[\begin{array}{ccccccc}
			\StMat_{\IntI\IntI} & \StMat_{\IntI \FSIIs}\\
			\StMat_{\FSIIs\IntI} & \StMat_{\FSIIs\FSIIs} &&&&&-(1-\tifs)\trans{\MoM}\\
			&&\FlMat_{\IntI\IntI}&\FlMat_{\IntI\FSIIs}&\gridsub{\FlMat}{\IntI\IntI}&\gridsub{\FlMat}{\IntI\FSIIs}\\
			&&\FlMat_{\FSIIs\IntI}&\FlMat_{\FSIIs\FSIIs}&\gridsub{\FlMat}{\FSIIs \IntI}&\gridsub{\FlMat}{\FSIIs\FSIIs}&(1-\tiff)\trans{\MoD}\\
			&&&&\AleMat_{\IntI\IntI}&\AleMat_{\IntI\FSIIs}\\
			&-\MoM&&\FlAle\MoD\\
      &\MoM&&&&-\MoD
		\end{array}\right]_{i}^{n+1}
\end{align}
and the residual vector
\begin{align}
\label{eq:rhsFSfull}
\rhsm_i^{n+1}
=	\left[\begin{array}{c}
		\strsub{\rhsm}{\IntI}\\
		\strsub{\rhsm}{\FSIIs}\\
		\flusub{\rhsm}{\IntI}\\
		\flusub{\rhsm}{\FSIIs}\\
		\gridsub{\rhsm}{\FSIIs}\\
    	\mao 0\\
		\mao 0
 	\end{array}\right]_i^{n+1}
+	\left[\begin{array}{c}
		\mao 0\\
		-\tifs\trans{\MoM}\mao\blamb^n\\
		\mao 0\\
		\tiff\trans{\MoD}\mao\blamb^n\\
		\mao 0\\
    \mao 0\\
		\mao 0
 	\end{array}\right]
+	\delta_{i0}
	\left[\begin{array}{c}
		\mao 0\\
		\mao 0\\
		\mao 0\\
		\mao 0\\
    \mao 0\\
		\Dt\MoD\,\flusupsub{\mao u}{n}{\FSIIs} - \MoM\strsub{\Delta \mao d}{\FSIIs,p}\\
    \MoM\strsub{\Delta \mao d}{\FSIIs,p}\\
	\end{array}\right].
\end{align}
\end{subequations}

The sixth row that brings the discrete kinematic constraint~\eqref{eq:DisVelConversion} into the system can be resolved for the fluid interface velocity increment
\begin{align}
\label{eq:KinemConstrFS}
  \flusupsub{\Delta \mao u}{n+1}{\FSIIs} 
  & = \FlAleI \MoP \strsupsub{\Delta \mao d}{n+1}{\FSIIs}
    + \delta_{i0}\FlAleI \MoP \strsub{\Delta \mao d}{\FSIIs,p}   
    - \delta_{i0} \frac{\Dt}{\FlAle}\flusupsub{\mao u}{n}{\FSIIs}
\end{align}
with the mortar projection matrix~$\MoP$ defined in~\eqref{eq:MortarP}. The coupling of structure displacements and ALE displacements at the fluid-structure interface is given by the last row in~\eqref{eq:SystemFSfull} and reads
\begin{align}
  \label{eq:ALEConstrFS}
  \Delta\gridsupsub{\mao d}{n+1}{\FSIIs,i+1}
  & = \MoP\Delta\strsupsub{\mao d}{n+1}{\FSIIs,i+1}
    + \delta_{i0}\,\MoP\Delta\strsub{\mao d}{\FSIIs,p}.
\end{align}
From the balance of linear momentum of the fluid interface degrees of freedom together with~\eqref{eq:KinemConstrFS} and~\eqref{eq:ALEConstrFS}, the unknown \Lagr multipliers are expressed by
\begin{align}
\label{eq:RecoverLambdaFS}
\begin{split}
	\blamb^{n+1} 
  &	= - \frac{\tiff}{1-\tiff}\mao\blamb^{n}
    - \frac{1}{1-\tiff}\itra{\MoD}
			\left(\flusupsub{\rhsm}{n+1}{\FSIIs}
          + \left(\FlAleI\FlMat_{\FSIIs\FSIIs}+\gridsub{\FlMat}{\FSIIs\FSIIs}\right)\MoP\Delta\strsupsub{\mao d}{n+1}{\FSIIs,i+1}\right)\\
  & - \frac{1}{1-\tiff}\itra{\MoD}
			\left(\FlMat_{\FSIIs\IntI}\Delta\flusupsub{\mao u}{n+1}{\IntI,i+1}
          + \gridsub{\FlMat}{\FSIIs\IntI}\Delta\gridsupsub{\mao d}{n+1}{\IntI,i+1}\right)\\
  & - \delta_{i0}\,\frac{1}{1-\tiff}\itra{\MoD}
      \left(\left(\FlAleI\FlMat_{\FSIIs\FSIIs}+\gridsub{\FlMat}{\FSIIs\FSIIs}\right)\MoP\Delta\strsub{\mao d}{\FSIIs,p} - \frac{\Dt}{\FlAle}\FlMat_{\FSIIs\FSIIs}\flusupsub{\mao u}{n}{\FSIIs}\right).
\end{split}
\end{align}
Equation~\eqref{eq:RecoverLambdaFS} is used to recover the \Lagr multiplier solution at the end of each time step as a postprocessing step.

Using~\eqref{eq:KinemConstrFS}, \eqref{eq:ALEConstrFS}, and~\eqref{eq:RecoverLambdaFS}, we are able to condense the saddle-point type system of equations. Finally, the condensed linear system for the structure-handled interface motion consists of the Jacobian matrix
\begin{subequations}
\label{eq:SystemFScond} 
\begin{align}
\label{eq:JacFScond}
\FSIMat^{n+1}_i
& = \left[\begin{array}{cccc}
			\StMat_{\IntI\IntI} & \StMat_{\IntI\FSIIs}\\
			\StMat_{\FSIIs\IntI} &  \StMat_{\FSIIs\FSIIs}+\frac{1-\tifs}{1-\tiff}\FlAleI\trans{\MoP} \FlMat_{\FSIIs\FSIIs}\MoP + \frac{1-\tifs}{1-\tiff}\trans{\MoP}\gridsub{\FlMat}{\FSIIs\FSIIs}\MoP & \frac{1-\tifs}{1-\tiff}\trans{\MoP}\FlMat_{\FSIIs \IntI} & \frac{1-\tifs}{1-\tiff}\trans{\MoP}\gridsub{\FlMat}{\FSIIs \IntI}\\
			&\FlAleI\FlMat_{\IntI\FSIIs}\MoP+\gridsub{\FlMat}{\IntI\FSIIs}\MoP & \FlMat_{\IntI\IntI} & \gridsub{\FlMat}{\IntI\IntI}\\
			&\AleMat_{\IntI\FSIIs}\MoP && \AleMat_{\IntI\IntI}
		\end{array}\right]_{i}^{n+1},
\end{align}
the solution increment vector
\begin{align}
\label{eq:SolFScond}
\trans{{\Delta\mao x_{i+1}^{n+1}}}
& =	\left[\begin{array}{ccccc}
			\trans{{\strsub{\Delta \mao d}{\IntI}}}&
			\trans{{\strsub{\Delta \mao d}{\FSIIs}}}&
			\trans{{\flusub{\Delta \mao u}{\IntI}}}&
			\trans{{\gridsub{\Delta \mao d}{\IntI}}}
		\end{array}\right]_{i+1}^{n+1},
\end{align}
and the residual vector
\begin{align}
\label{eq:rhsFScond}
\begin{split}
\rhsm_i^{n+1}
&	=	\left[\begin{array}{c}
			\strsub{\rhsm}{\IntI}\\
			\strsub{\rhsm}{\FSIIs} + \frac{1-\tifs}{1-\tiff}\trans{\MoP}\flusub{\rhsm}{\FSIIs}\\
			\flusub{\rhsm}{\IntI}\\
			\gridsub{\rhsm}{\FSIIs}
		\end{array}\right]_i^{n+1}
	+	\left[\begin{array}{c}
			\mao 0\\
			\left(-\tifs+\frac{\tiff(1-\tifs)}{1-\tiff}\right)\trans{\MoM}\mao\blamb^n\\
			\mao 0\\			
			\mao 0
		\end{array}\right]
\\
&	+	\delta_{i0}
		\left[\begin{array}{c}
			\mao 0\\
			\frac{1-\tifs}{1-\tiff}\FlAleI\trans{\MoP}\FlMat_{\FSIIs\FSIIs}\MoP\Delta\strsub{\mao d}{\FSIIs,p} + \frac{1-\tifs}{1-\tiff}\trans{\MoP}\gridsub{\FlMat}{\FSIIs\FSIIs}\MoP\Delta\strsub{\mao d}{\FSIIs,p} - \frac{1-\tifs}{1-\tiff}\frac{\Dt}{\FlAle}\trans{\MoP}\FlMat_{\FSIIs\FSIIs}\flusupsub{\mao u}{n}{\FSIIs}\\
      \FlAleI\FlMat_{\IntI\FSIIs}\MoP\Delta\strsub{\mao d}{\FSIIs,p} + \gridsub{\FlMat}{\IntI\FSIIs}\MoP\Delta\strsub{\mao d}{\FSIIs,p} - \frac{\Dt}{\FlAle}\FlMat_{\IntI\FSIIs}\flusupsub{\mao u}{n}{\FSIIs}\\
      \AleMat_{\IntI\FSIIs}\MoP\Delta\strsub{\mao d}{\FSIIs,p}
		\end{array}\right].
\end{split}
\end{align}
\end{subequations}

\begin{remark}
\label{rem:InterfaceDBCsFS}
  As already indicated in remark~\ref{rem:InterfaceDBCsSS}, one has to be careful with essential boundary conditions at the fluid-structure interface~$\FSII$. Following the same arguments as before, now, only the structure side of the interface is allowed to carry Dirichlet boundary conditions. They will be imposed on the fluid side weakly via the mortar coupling.
\end{remark}

\section{Numerical examples}
\label{sec:examples}

Three numerical examples are used to demonstrate and discuss some properties of the presented solution schemes. First, a very simple test case is proposed, where an analytical solution is known, and used to study temporal convergence as well as some aspects of Dirichlet boundary conditions at the fluid-structure interface in~\S\ref{ssec:fp}. In~\S\ref{ssec:dc}, the well-known 2D driven cavity with flexible bottom is utilized to demonstrate the effect of predictors on the overall computational costs. Finally, different combinations of time integration schemes are compared to each other using the 3D pressure wave example mimicking hemodynamic conditions in~\S\ref{ssec:pw}.

In all three examples, equal-order interpolated linear finite elements with residual-based stabilization are used for spatial discretization of the fluid field. The structure field is discretized with mixed/hybrid finite elements. Enhanced assumed strains (EAS) are utilized to deal with locking phenomena.

\subsection{Pseudo 1D FSI}
\label{ssec:fp}

\begin{figure}
\begin{center}
  \subfigure[problem setup]{\label{fig:fp}\includegraphics[width=0.35\textwidth]{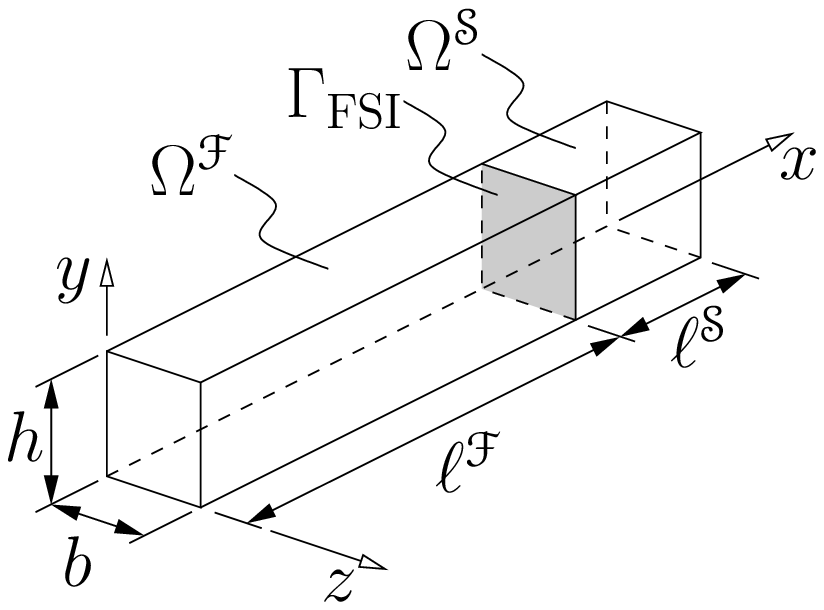}}
  \hfill
  \subfigure[solution of~$\flu{p}$ and~$\blamb$]{\label{fig:fp_lambda}\includegraphics[width=0.56\textwidth]{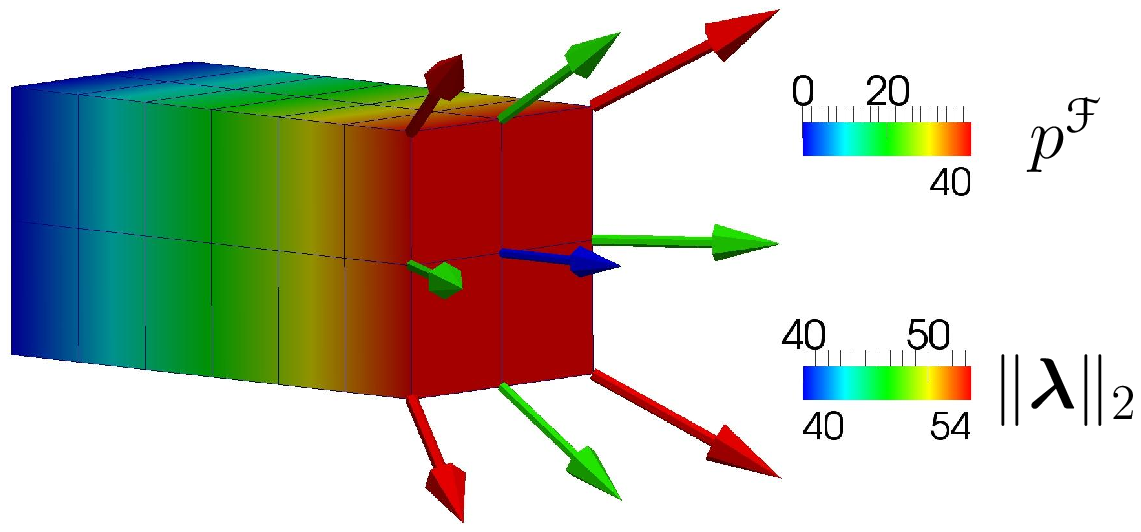}}
\caption{Geometry an solution of pseudo 1D FSI example with analytical solution --- \emph{Left:} The structural block~$\str{\Omega}$ moves in $x$-direction due to a time dependent Dirichlet boundary condition at~$x=\flu{\ell}+\str{\ell}$. Thus, fluid is pushed out or sucked in across the Neumann boundary at~$x=0$. All movement in $y$- and $z$-direction is suppressed, leaving a pseudo 1D problem. \emph{Right:} Pressure field in fluid domain~$\flu{\Omega}$ and \Lagr multiplier field~$\mao{\blamb}$: The \Lagr multiplier field represents the interface traction onto the structure. The $x$-components represent the fluid pressure exerted onto the structure. The lateral components in $y$- and $z$-direction constrain the $y$- and $z$-components of the fluid velocity.}
\end{center}
\vspace{-5mm}
\end{figure}

As a very simple example, we first consider a pseudo one-dimensional FSI problem as sketched in figure~\ref{fig:fp}.
It is used to demonstrate temporal convergence properties of the proposed monolithic solution scheme employing the comparison to an analytical solution. Additionally, the special role of Dirichlet boundary conditions at the fluid-structure interface as discussed in remarks~\ref{rem:InterfaceDBCsSS} and~\ref{rem:InterfaceDBCsFS} is illustrated by a visualization of the \Lagr multiplier field.

The example is set up as a real 3D problem, but is constrained to one dimension via Dirichlet boundary conditions, \ie all displacement and velocity degrees of freedom are forced to zero in $y$- and $z$-direction. Hence, movement is possible only in $x$-direction. The problem is driven by a time dependent Dirichlet boundary condition on the dry side of the solid domain~$\str{\Omega}$, \ie at~$x=\flu{\ell}+\str{\ell}$. When moving the structural block~$\str{\Omega}$, the size of the fluid volume~$\flu{\Omega}$ changes and fluid is pushed out or sucked in across the fluid Neumann boundary at~$x=0$. Assuming a time dependent Dirichlet boundary condition~$\strsub{\bar{\teo d}}{\mr{D}}(t)$ for the displacement of the dry side of the structure at~$x=\flu{\ell}+\str{\ell}$ as well as incompressible solid and fluid domains, the analytical solution for velocity field~$\teo u(x,t)$, acceleration field~$\teo a(x,t)$ and fluid pressure field~$\flu{p}(x,t)$ reads:
\label{eq:AnalyticalSolution}
\begin{align}
  \teo u(x,t) = \pDer{\strsub{\bar{\teo d}}{\mr{D}}}{t},\quad
  \teo a(x,t) = \spDer{\strsub{\bar{\teo d}}{\mr{D}}}{t},\quad
  \flu{p}(x,t) = - \flu{\rho}\flusub{\teo a}{x}\cdot x + \left.\flusub{p}{\infty}\right\vert_{x=0}
\end{align}

In this example, the structure field is chosen as master field (\cf~\S\ref{ssec:FluidSplit}). Hence, at the interface, only the structural degrees of freedom are subject to Dirichlet boundary conditions. According to remark~\ref{rem:InterfaceDBCsFS}, the fluid side is not allowed to carry Dirichlet boundary conditions.

Due to the spatially constant velocity field and the spatially linear pressure field, the finite element solution can capture the spatial distribution of the analytical solution exactly.

When choosing the imposed time dependent Dirichlet boundary condition~$\strsub{\bar{\teo d}}{\mr{D}}(t)$ such that the analytical solution is also contained in the discrete temporal solution space, for example~$\strsub{\bar{\teo d}}{\mr{D}}(t) = -t^2$, the analytical solution is fully recovered by the numerical scheme up to machine precision.

\begin{figure}
\begin{center}
  \subfigure[$\ltwo$-error in velocity field]{\label{fig:fp_velocity}\includegraphics[width=0.48\textwidth]{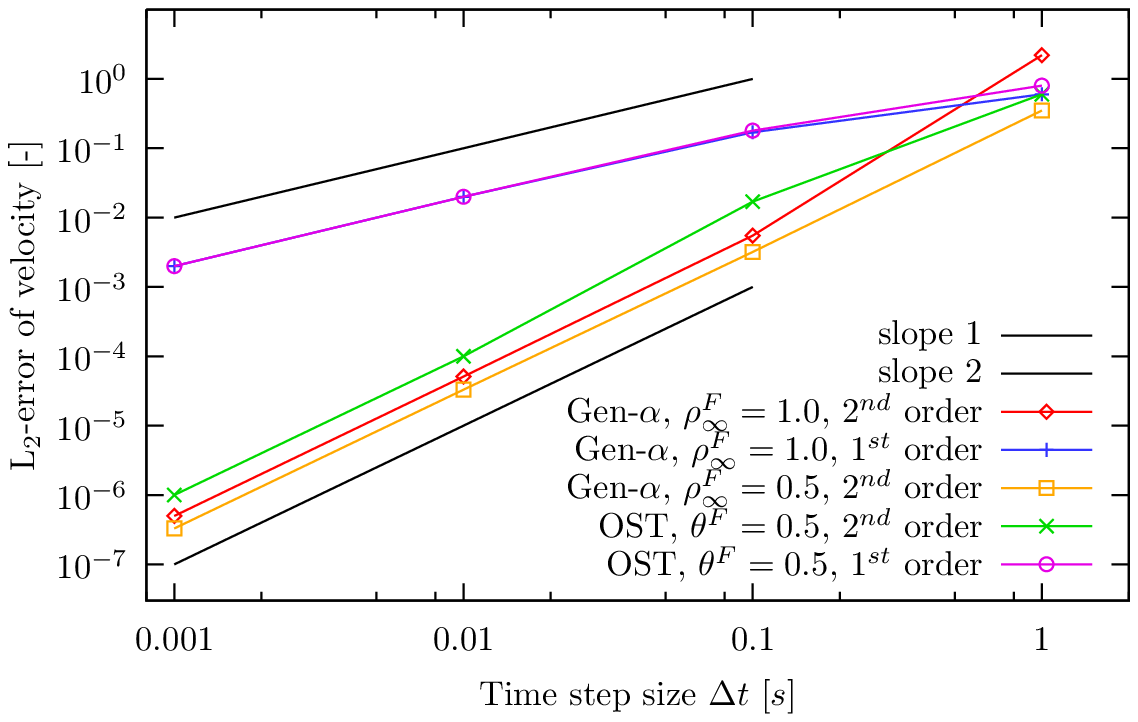}}
  \hfill
  \subfigure[$\ltwo$-error in pressure field]{\label{fig:fp_pressure}\includegraphics[width=0.48\textwidth]{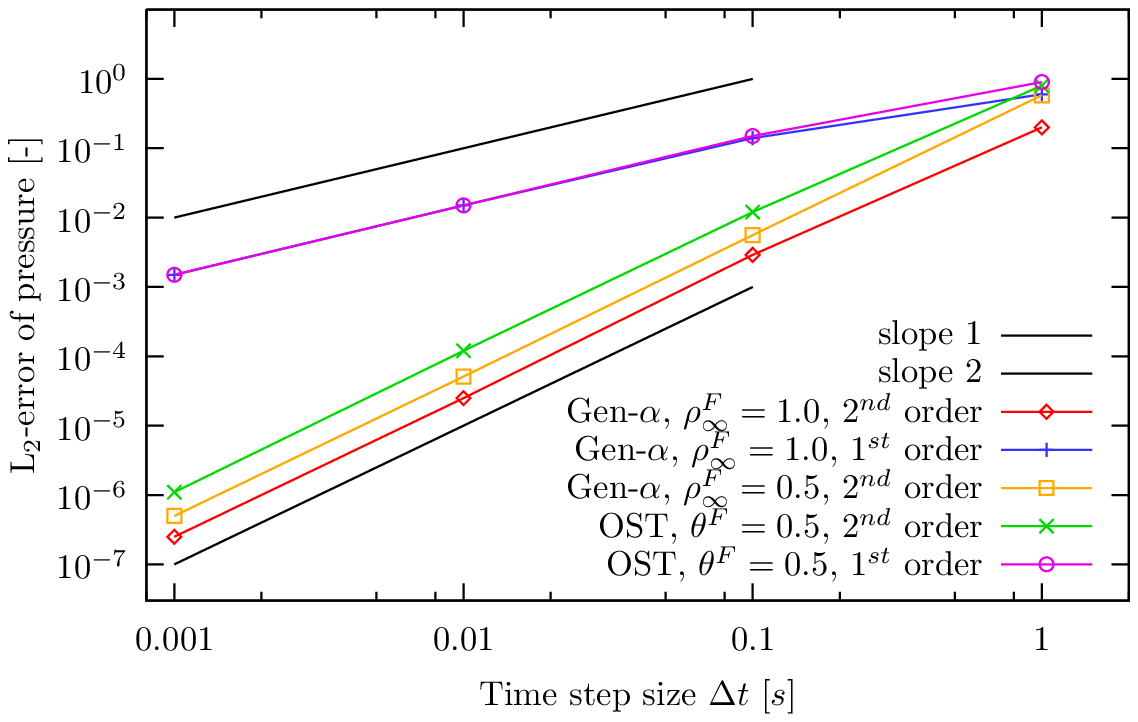}}
\caption{Temporal convergence study for pseudo 1D FSI example --- Comparison of temporal convergence for different fluid time integration schemes and different conversions of ALE interface displacements into fluid interface velocities. Temporal convergence is measured using $\ltwo$-errors of velocity field and pressure field in the fluid volume~$\flu{\Omega}$. The computed convergence rates match theoretical expectations perfectly.}
\label{fig:fp_timeconvergence}
\end{center}
\vspace{-2mm}
\end{figure}

For this reason, in order to study temporal convergence, a Dirichlet boundary condition~$\strsub{\bar{\teo d}}{\mr{D}}(t) = -t^5$ is prescribed on the structure. The spatial solution can still be fully recovered, but the involved time integration schemes are not able to capture the temporal evolution exactly. Hence, temporal refinement should lead to error reduction. For this study, we calculate the $\ltwo$-error of the velocity and pressure field in the fluid volume~$\flu{\Omega}$ compared to the analytical solution~\eqref{eq:AnalyticalSolution}. The actual material parameters are of no importance. For temporal discretization of the structure field, generalized-$\alpha$ time integration with spectral radius~$\strsub{\rho}{\infty}=1.0$, \ie without numerical dissipation, is used. The fluid time integrator is either the generalized-$\alpha$ scheme with various spectral radii~$\flusub{\rho}{\infty}$ or the one-step-$\theta$ scheme with various choices for~$\flu{\theta}$. The conversion between ALE displacements and fluid velocities is varied between trapezoidal rule and backward Euler as indicated in~\eqref{eq:FluAleConversion} and~\eqref{eq:DefinitionTau}. When structure and fluid time integration scheme as well as the conversion between ALE displacements and fluid velocities are chosen to be second order accurate, the overall FSI scheme is expected to be second order accurate in time as well. As soon as one of them is only first order accurate in time, the order of temporal accuracy of the overall algorithm is expected to reduce to first order. Figure~\ref{fig:fp_timeconvergence} shows the temporal convergence plots for velocity field and pressure field in the fluid volume~$\flu{\Omega}$. Time step sizes and fluid time integration schemes with particular parameters are detailed in figure~\ref{fig:fp_timeconvergence}.
Only the cases where the conversion of ALE displacements and fluid velocities is done with the backward Euler formula~\eqref{eq:backwardeuler} the temporal convergence order deteriorates to first order. In all other cases, where the overall algorithm is expected to be second order accurate, the temporal convergence is of second order. Altogether, the theoretically expected convergence rates are fully recovered by the proposed monolithic FSI scheme.

Looking at the Dirichlet boundary conditions at the fluid-structure interface, that have to preclude lateral motions, one has to consider that interface degrees of freedom on the slave side, \ie in the fluid field, are not allowed to carry Dirichlet boundary conditions (\cf remark~\ref{rem:InterfaceDBCsFS}). Thus, the fluid interface degrees of freedom are not constrained by Dirichlet boundary conditions at all. As such, the Dirichlet boundary conditions of the structural interface degrees of freedom are assigned to the fluid interface degrees of freedom weakly via the mortar coupling. The traction that forces the $y$- and $z$-components of the fluid velocity to zero is represented by the $y$- and $z$-components of the \Lagr multiplier field. Figure~\ref{fig:fp_lambda} shows a visualization of the fluid domain and the \Lagr multiplier field. The weak enforcement of Dirichlet boundary conditions across the fluid-structure interface is clearly observed, since the traction field exhibits components in lateral $y$- and $z$-direction.

\subsection{Driven cavity with flexible bottom}
\label{ssec:dc}

\begin{figure}
\begin{center}
  \subfigure[geometry and boundary conditions]{\label{fig:dc}\includegraphics[width=0.42\textwidth]{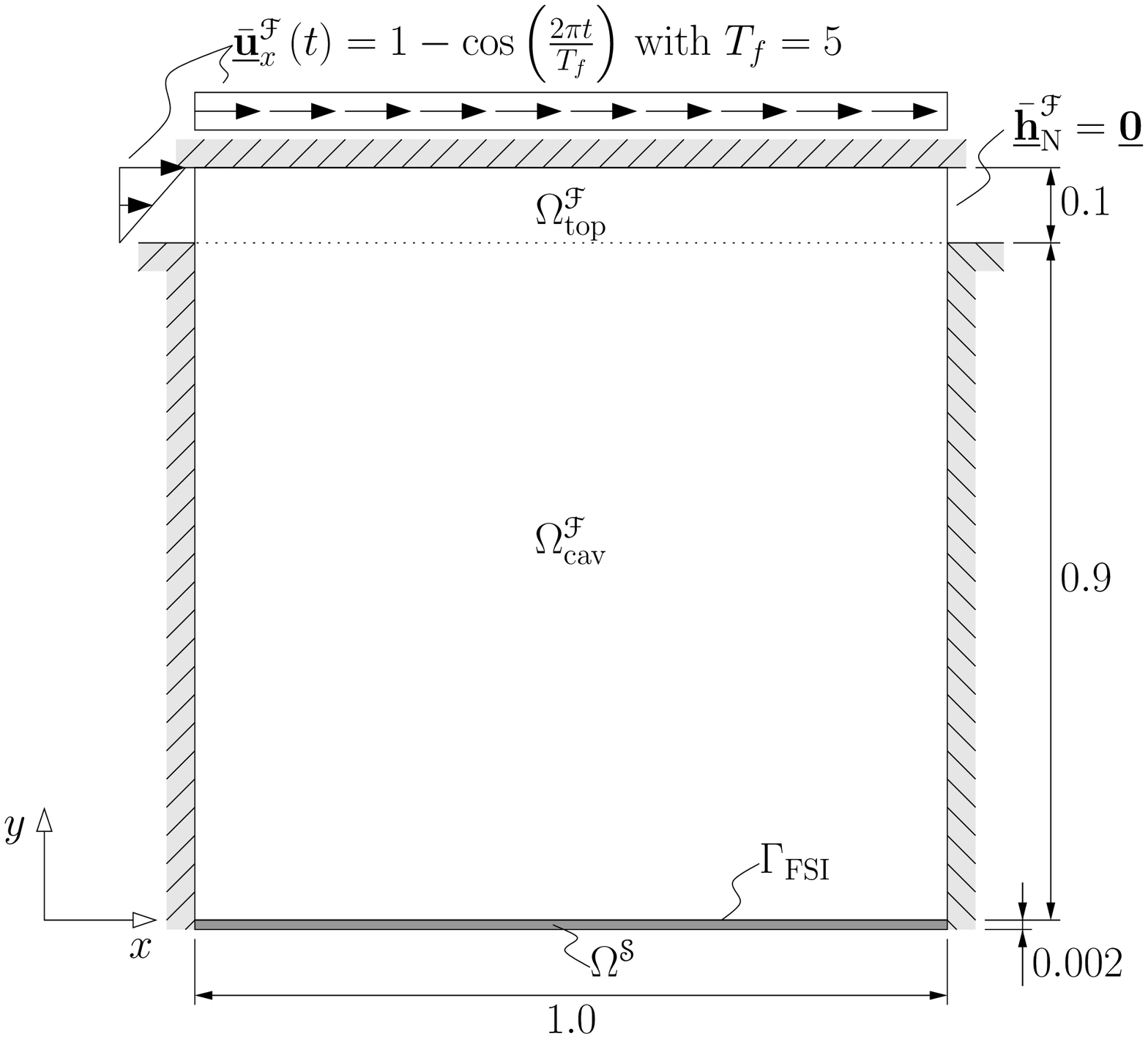}}
  \hfill
  \subfigure[solution at time~$t=19.0$]{\label{fig:dc_solution}\includegraphics[width=0.53\textwidth]{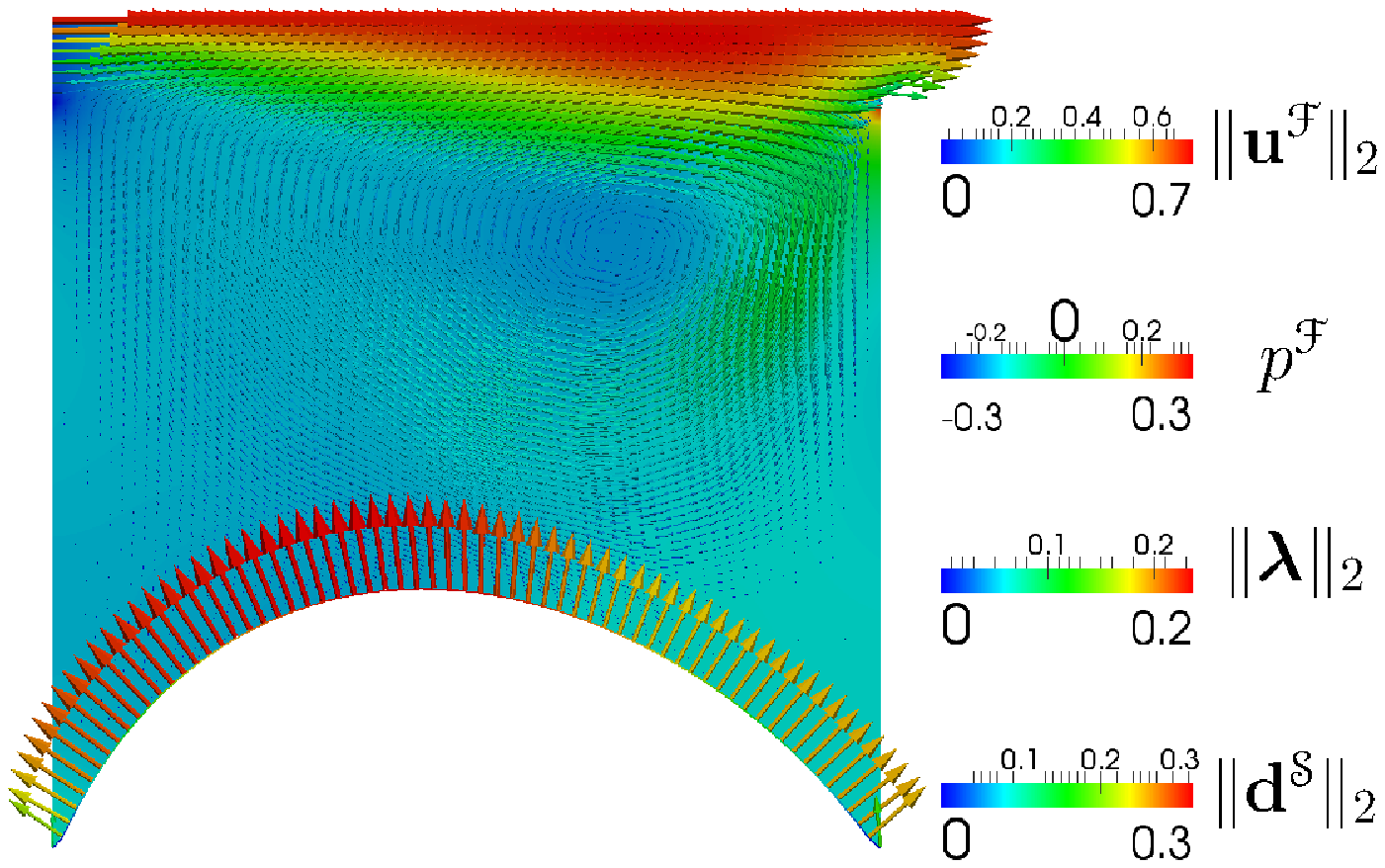}}
\caption{Geometry, boundary conditions, and solution for leaky driven cavity with flexible bottom --- \emph{Left:} The square fluid domain~$\flu{\Omega}$ is split into a cavity portion~$\flusub{\Omega}{\mr{cav}}$, covering the bottom part of~$\flu{\Omega}$, and a top portion~$\flusub{\Omega}{\mr{top}}$. On left and right walls of the cavity, \emph{no-slip} boundary conditions are imposed. On the top of~$\flusub{\Omega}{\mr{top}}$, the velocity in $x$-direction is prescribed by~$\flusub{\bar{\teo u}}{x}\left(t\right)$, whereas the velocity in $y$-direction is set to zero. The left side of the top region is subject to a linearily varying prescribed inflow velocity. The pressure level is determined by a do-nothing Neumann boundary condition on the right side of the top region. The structural domain~$\str{\Omega}$ is clamped on its left and right edges.  \emph{Right:} The fluid domain is shown via a contour plot of the pressure field. Additionally, the fluid velocity field is visualized using a vector plot. At the interface, the \Lagr multiplier field, \ie the coupling traction, is shown as traction vectors.}
\end{center}
\vspace{-5mm}
\end{figure}

To demonstrate the effect of the predictors, a two-dimensional leaky driven cavity with flexible bottom as sketched in figure~\ref{fig:dc} is used (see also~\cite{Mok2001a}).
The structure is modelled with a St.-Venant-Kirchhoff material with Young's modulus~$\str{E}=250$, Poisson's ratio~$\str{\nu}=0$, and density~$\str{\rho}=500$. The incompressible Newtonian fluid has a dynamic viscosity of~$\flu{\viscdyn}=0.01$ and a density~$\flu{\rho}=1$. Geometry, dimensions, and boundary conditions are detailed in figure~\ref{fig:dc}.

Spatial discretization is performed with two different grids. For the coarser grid, the cavity volume~$\flusub{\Omega}{cav}$ is meshed with $64\times64$ bilinear quadrilateral elements, whereas the top volume~$\flusub{\Omega}{top}$ is meshed with $64\times8$ bilinear quadrilateral elements. In order to realize non-matching grids at the fluid-structure interface~$\FSII$, the structure is discretized with $72\times2$ bilinear quadrilateral elements. For the finer grid, the number of elements in each direction is doubled. For temporal discretization both fields employ generalized-$\alpha$ time integration without numerical dissipation and a time step size~$\Dt = 0.01$. Figure~\ref{fig:dc_solution} depicts the solution at time~$t=19.0$.

This example was ran using different types of predictors in the structure field only to demonstrate their effect on computational costs. The fluid field is always treated without any predictor. The reference solution is computed without any predictors, \ie assuming constant displacements, velocities, and accelerations (referred to as \emph{ConstDis}). In the structure field, two types of predictors are used: First, a constant structural velocity is assumed, yielding a linear displacement prediction (referred to as \emph{ConstVel}). Secondly, the accelerations are assumed to be constant, resulting in a linear extrapolated velocity field and a quadratically extrapolated displacement field (referred to as \emph{ConstAcc}). All predictors require only simple and extremely cheap vector operations like multiplication with a scalar and addition of vectors. Hence, they are of negligible costs compared to the remaining operations, especially to the costs of the linear solver.

The costs are quantified by the number of linear solver iterations per time step, since in general these costs dominate, especially when it comes to large problem sizes. Hence, a reduction of the number of linear solver iterations has huge impact on the overall computational costs and, thus, is very desirable.

For comparison, the simulation parameters except for the mesh size have been held constant. The linear solver utilizes an ILU(0) preconditioner for each field and is solved using a GMRES procedure~\cite{Saad1986} where the Krylov space dimension is set to~$50$. The iterative linear solver stops, when the relative residual~$||\rhsm||_2 / ||\rhsm_0||_2$ is below~$10^{-5}$. The nonlinear iteration is stopped, as soon as the residuals as well as the nonlinear solution increments of the displacement field, the velocity field, and the pressure field measured in~$\ltwo$- and~$\linf$-norm are below~$10^{-8}$. In addition, the $\ltwo$- and~$\linf$-norm of the interface residual and increment were required to be smaller than~$10^{-9}$.

\begin{figure}
\begin{center}
  \subfigure[$64\times 64$ mesh: reduction by 9.7\% on average]{\includegraphics[width=0.48\textwidth]{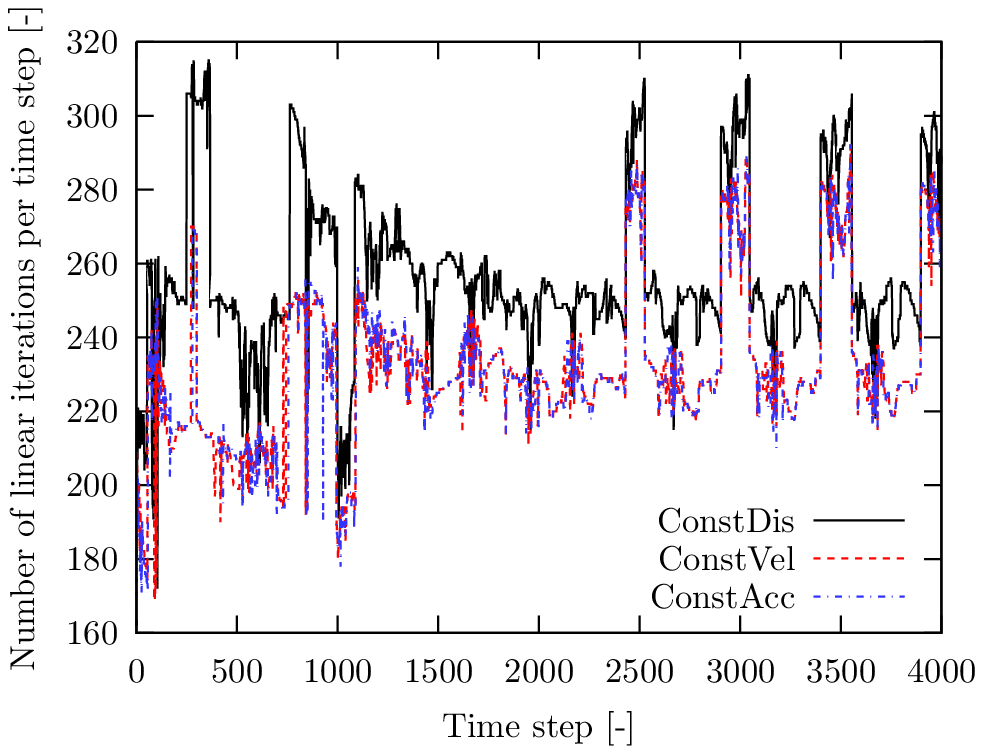}}\hfill
  \subfigure[$128\times 128$ mesh: reduction by 3.1\% on average]{\includegraphics[width=0.48\textwidth]{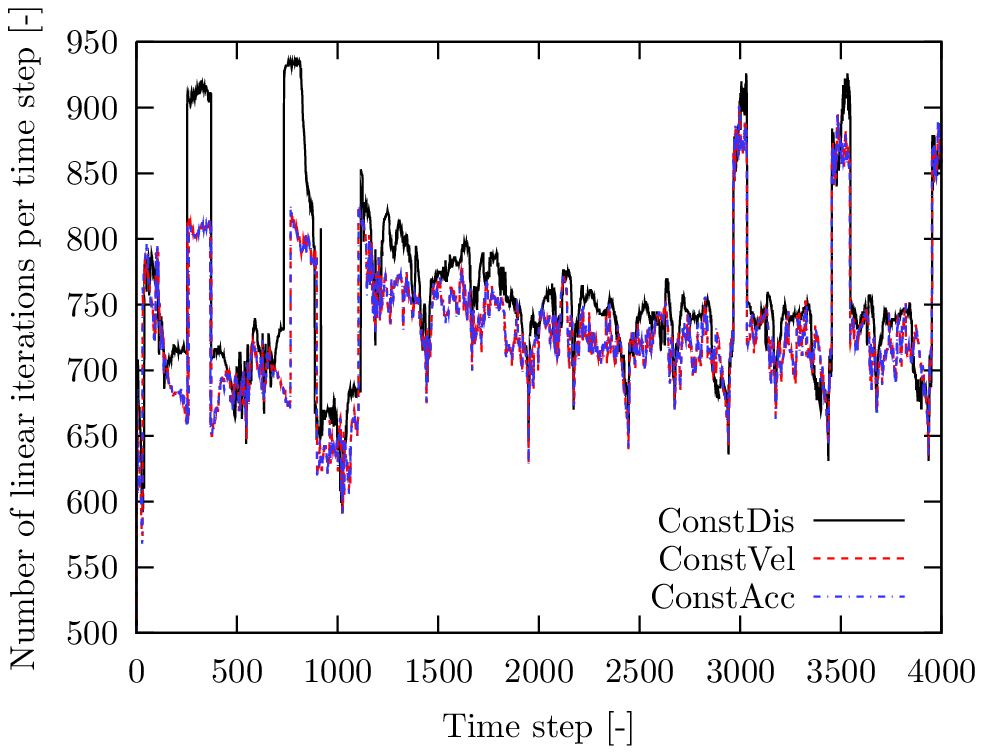}}
\caption{Driven cavity with flexible bottom --- Comparison of number of linear iterations per time step for different choices of predictors. Using simple predictors like \emph{ConstVel} or \emph{ConstAcc} in the structure field reduces the number of linear iterations per time step and, thus, the actual computational costs compared to the case without predictor, \ie \emph{ConstDis}.}
\label{fig:dc_liniters}
\end{center}
\vspace{-2mm}
\end{figure}

The number of linear iterations per time step for different choices of predictors is shown in figure~\ref{fig:dc_liniters}.
Just by employing the \emph{ConstVel} or the \emph{ConstAcc} predictor in the structure field, the average number of linear iterations per time step is reduced by 9.7\% on the coarser grid. On the fine grid, the reduction is 3.1\%.

\begin{remark}
In principle, the solution schemes proposed in~\S\ref{ssec:StruSplit} and~\S\ref{ssec:FluidSplit} are able to handle fluid predictors as well. However, if only the velocity field is predicted in a comparably simple way as in the structure field, the pressure field does not match the velocity field after the prediction. Hence, we recommend to just predict the structural solution unless sophisticated fluid predictors that include a pressure projection step are available. In our implementation, however, we only consider explicit predictors that come with only negligible additional cost and therefore refrain from using a non-constant fluid predictor step here.
\end{remark}

When looking at the interface energy production per step as discussed in remark~\ref{rem:InterfaceEnergy}, the amount of energy per step~$\Delta E_{\FSIIs}^{n\rightarrow n+1}$ in the worst case scenario is at the order of~$10^{-6}$ whereas the kinetic energy of the system is at the order of~$10^{-1}$. Furthermore we keep in mind the physical dissipation due to fluid viscosity and further numerical dissipation stemming from time integration or fluid stabilization. Hence, the energy production per step can be considered as negligible. For~$\tifs=\tiff$, $\Delta E_{\FSIIs}^{n\rightarrow n+1}$ vanishes up to machine precision. For this reason, the observations in remark~\ref{rem:InterfaceEnergy} are confirmed numerically.

\subsection{Pressure wave through collapsible tube}
\label{ssec:pw}

Finally, a pressure wave travelling through a collapsible tube (see \eg~\cite{Gee2011,Gerbeau2003}) is examined mimicking hemodynamic conditions. The outstanding efficiency of a monolithic solution scheme compared to partitioned schemes has already been demonstrated in~\cite{Kuettler2010}. A more detailed analysis of the performance of the linear solvers has been performed in~\cite{Gee2011} for this example. 

\begin{figure}
\begin{center}
  \includegraphics[width=0.4\textwidth]{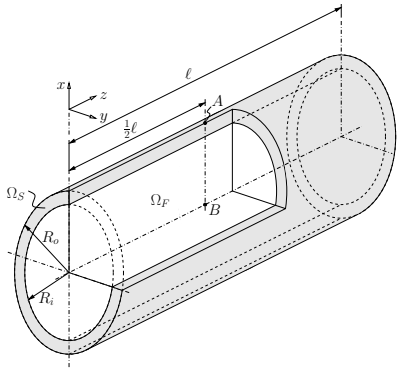}
\caption{Geometry of pressure wave example --- A solid tube (Young's modulus~$\unit{E=3\cdot 10^6}{g/(cm\cdot s^2)})$, Poisson's ratio~$\str{\nu}=0.3$, density~$\str{\rho} = \unit{1.2}{g/cm^3}$, outer radius~$R_o=\unit{0.6}{cm}$, inner radius~$R_i=\unit{0.5}{cm}$, length~$\ell=\unit{5.0}{cm}$) is filled with an incompressible Newtonian fluid (dynamic viscosity~$\flu{\viscdyn}=\unit{0.03}{g/(cm\cdot s)}$, density~$\flu{\rho}=\unit{1.0}{g/cm^3}$) that is initially at rest.}
\label{fig:pw}
\end{center}
\vspace{-5mm}
\end{figure}

The geometry is depicted in figure~\ref{fig:pw}.
The solid tube is clamped at both ends. The fluid is initially at rest. For the duration of~$\unit{3\cdot 10^{-3}}{s}$, it is loaded with a surface traction~$\flu{\bar{\teo h}}=\unit{1.3332\cdot 10^4}{g\cdot cm/s^2}$ in $z$-direction at~$z=0$. At~$z=\ell$, fluid velocities are prescribed to zero, meaning that the tube is closed at that end. As a result, a pressure wave travels along the tube's longitudinal axis and is reflected at the closed end of the tube. The constitutive behavior of the structure is modeled by a St.-Venant-Kirchhoff material, the fluid is assumed to be an incompressible Newtonian fluid. The actual material parameters are given in figure~\ref{fig:pw}.

The solid is discretized with trilinear hexahedral elements and utilizes \emph{enhanced assumed strains (EAS)} in order to circumvent possible locking phenomena. For the fluid discretization, trilinear stabilized equal-order hexahedral elements are used. The problem is solved on five different grids \emph{pw1}-\emph{pw5} with different levels of mesh refinement. 
Table~\ref{tab:pw_dofs} provides the number of unknowns per field for each of the five meshes.
\begin{table}
\footnotesize
\caption{Pressure wave trough collapsible tube --- number of unknowns for five different meshes \emph{pw1}-\emph{pw5}}
\label{tab:pw_dofs}
\begin{center}
\begin{tabular}{c|ccccccc}
  \hline
  Mesh & \# of structure DOFs & \# of fluid DOFs & \# of ALE DOFs & Total \# of DOFs\\
  \hline
  \emph{pw1} & $2016$ & $7476$ & $5607$ & $15099$\\
  \emph{pw2} & $11808$ & $55268$ & $41451$ & $108527$\\
  \emph{pw3} & $23424$ & $181780$ & $136335$ & $341539$\\
  \emph{pw4} & $77760$ & $425412$ & $319059$ & $822231$\\
  \emph{pw5} & $556416$ & $3339140$ & $2504355$ & $6399911$\\
  \hline
\end{tabular}
\end{center}
\end{table}
Temporal discretization is performed with generalized-$\alpha$ time integration in the structure field and generalized-$\alpha$ or one-step-$\theta$ time integration in the fluid field.  In all cases, the displacement-velocity conversion at the fluid-structure interface is done with the second order accurate trapezoidal rule~\eqref{eq:trapezoidal}. 

Mesh independence is examined. In addition, the effect of different combinations of time integration schemes and time integration parameters~$\strsub{\rho}{\infty}$,~$\flusub{\rho}{\infty}$, and~$\flu{\theta}$ in both structure and fluid field is studied in detail. Furthermore, solutions obtained with different time step sizes~$\Dt$ are compared to each other. For comparison, we monitor the temporal evolution of the radial displacement~$\strsub{\ensuremath{d}}{x}$ at point~$A(\unit{0.6}{cm},0,\unit{2.5}{cm})$ on the one hand. On the other hand, the temporal evolution of the fluid pressure~$\flu{p}$ at the center point~$B(0,0,\unit{2.5}{cm})$ is observed.

\begin{figure}
\begin{center}
  \subfigure[radial displacement]{\includegraphics[width=0.48\textwidth]{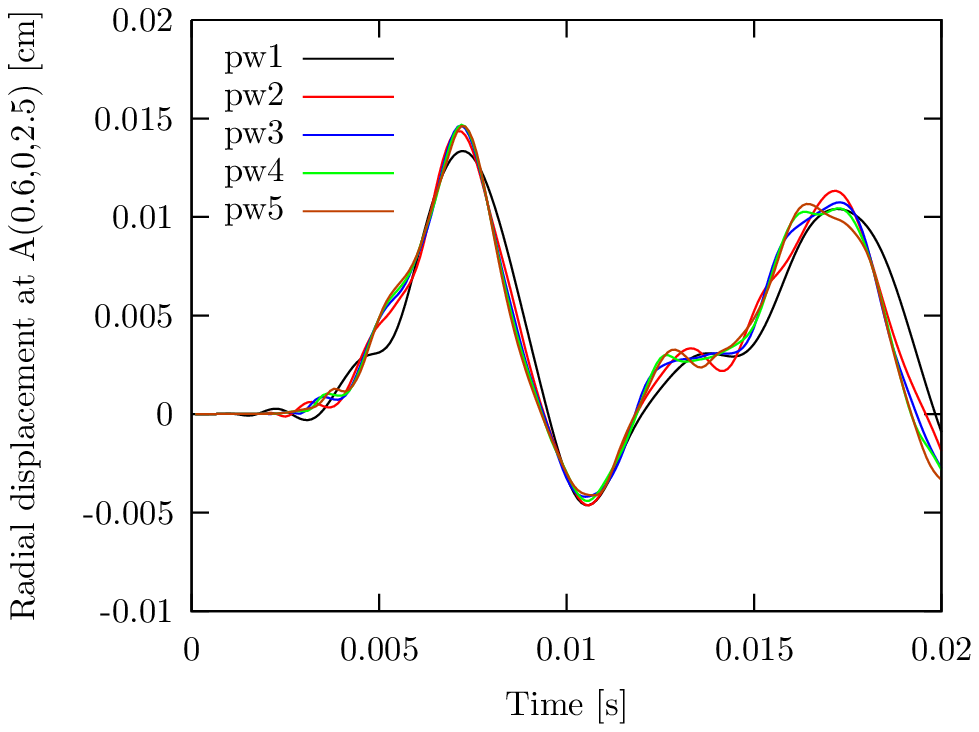}}
  \hfill
  \subfigure[fluid pressure]{\includegraphics[width=0.48\textwidth]{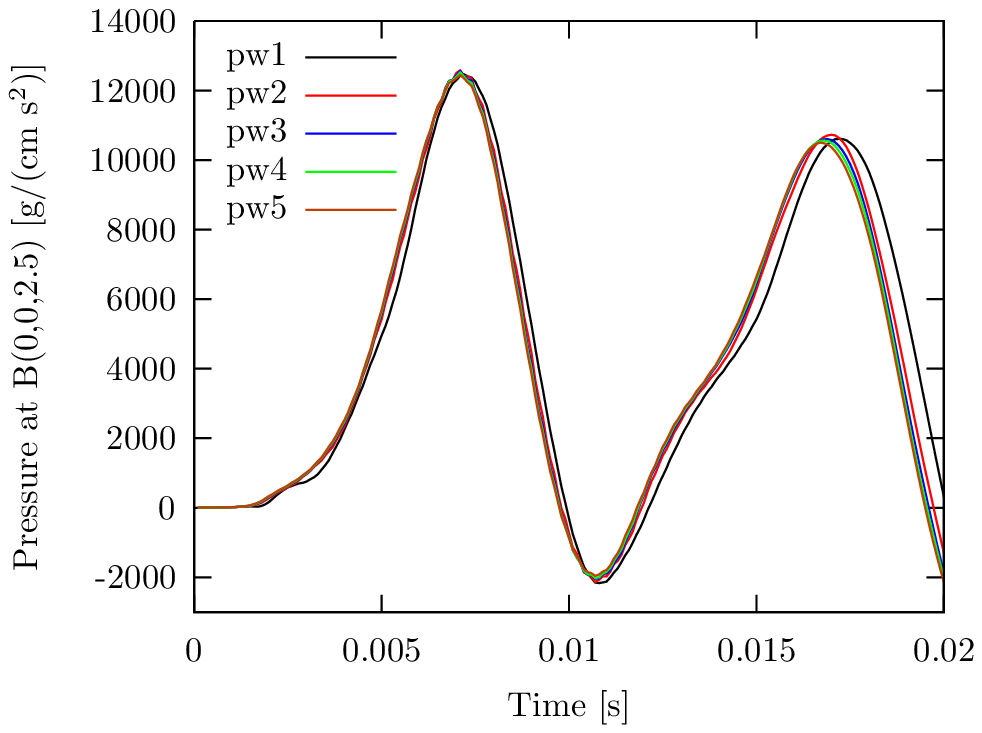}}
\caption{Mesh refinement for pressure wave example --- It can be seen that the influence of the mesh on the overall behavior of the solution is small for the meshes \emph{pw2}-\emph{pw5}. Only the solution obtained on the coarsest grid~\emph{pw1} differs significantly.}
\label{fig:pw_h}  
\end{center}
\vspace{-5mm}
\end{figure}

In figure~\ref{fig:pw_h}, the solutions for the different meshes \emph{pw1}-\emph{pw5} are reported.
Therefore, all computations have been carried out with a time step size~$\Dt=\unit{1.0\cdot10^{-4}}{s}$ as it is usual in literature~\cite{Gee2011,Kloeppel2011,Kuettler2010} and generalized-$\alpha$ time integration in both fields. The spectral radii have been chosen to~$\strsub{\rho}{\infty}=0.8$ and~$\flusub{\rho}{\infty}=0.5$ for the structure and fluid time integrator, respectively.
The first peak occurs, when the pressure wave passes points~$A$ and~$B$ for the first time. After the reflection at the closed end of the tube, the pressure wave travels in negative $z$-direction and causes the second peaks in figure~\ref{fig:pw_h}. Only the solution on the coarsest grid~\emph{pw1} differs significantly from the other fine grid solutions. Thus, for all further comparisons, we use the medium-sized discretization~\emph{pw3}. 

\begin{figure}
\begin{center}
  \subfigure[radial displacement]{\includegraphics[width=0.48\textwidth]{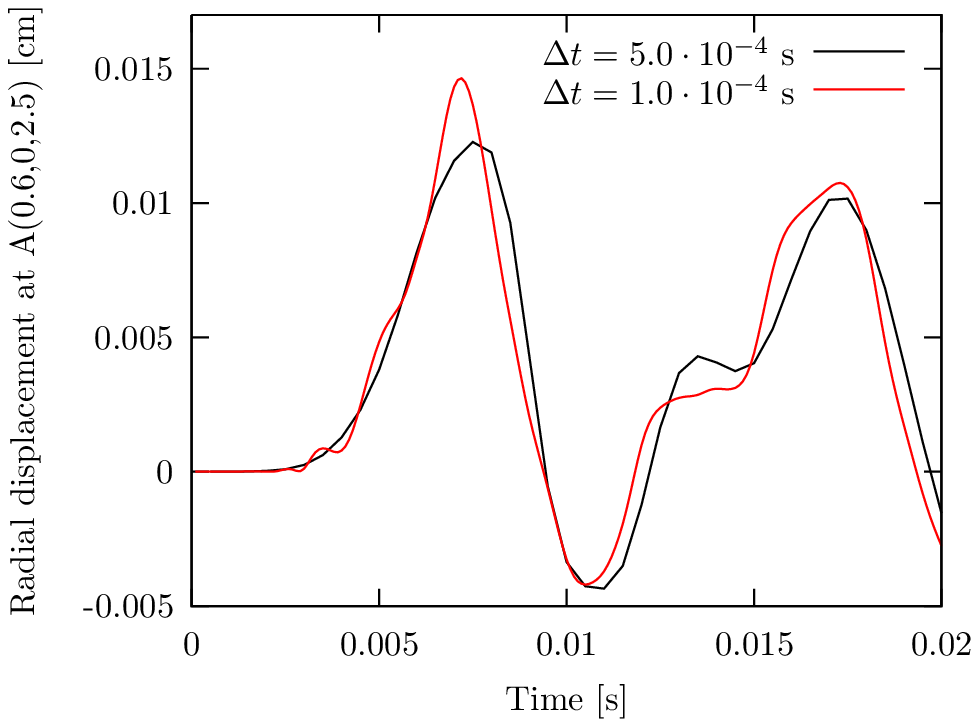}}
  \hfill
  \subfigure[fluid pressure]{\includegraphics[width=0.48\textwidth]{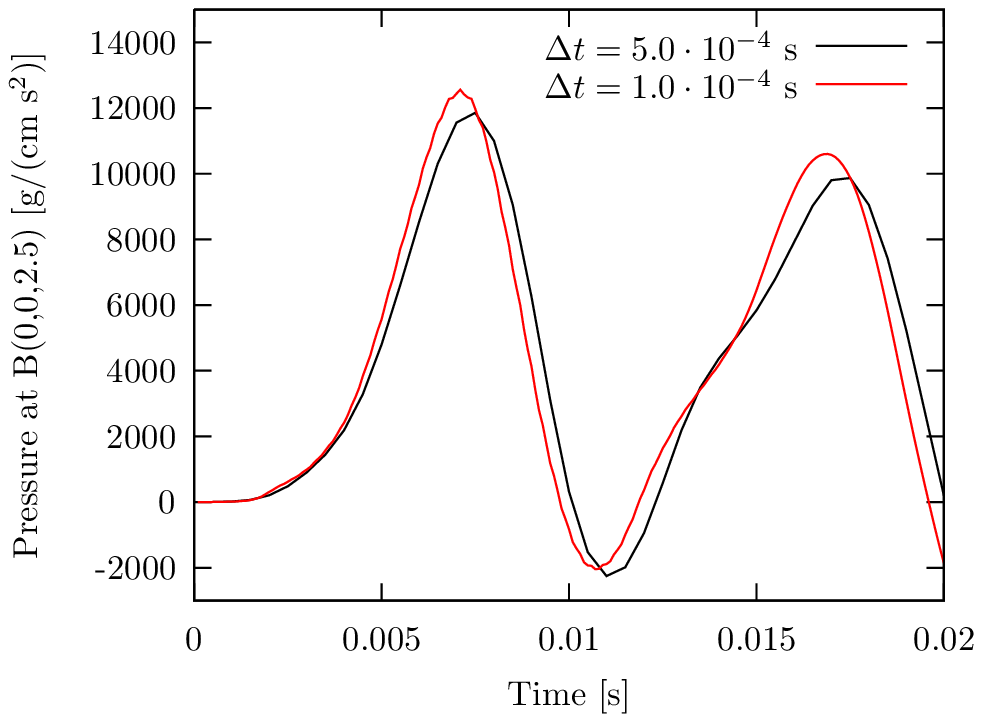}}
\caption{Different time step sizes for pressure wave example --- Using~$\Dt=\unit{1.0\cdot10^{-4}}{s}$ as reference solution, a larger time step size does not resolve the solution properly. }
\label{fig:pw_t}  
\end{center}
\vspace{-5mm}
\end{figure}

A comparison of different time step sizes has been performed and is reported in figure~\ref{fig:pw_t}.
Again, both fields have been integrated with generalized-$\alpha$ time integration with spectral radii~$\strsub{\rho}{\infty}=0.8$ and~$\flusub{\rho}{\infty}=0.5$ for the structure and fluid time integrator, respectively. The larger time step~$\Dt=\unit{5.0\cdot10^{-4}}{s}$ is not able to resolve the problem properly, leading to a shift of the displacement and pressure maxima to larger time values. The better solution is obtained with~$\Dt=\unit{1.0\cdot10^{-4}}{s}$. Due to the small time step limit of the fluid stabilization~\cite{Bochev2004}, a further refinement of the time step size is not possible without loosing the stabilizing effects of the fluid stabilization. In accordance with the literature~\cite{Gee2011,Kloeppel2011,Kuettler2010}, all further computations are done with a time step size of~$\Dt=\unit{1.0\cdot10^{-4}}{s}$.

\begin{figure}
\begin{center}
  \subfigure[radial displacement]{\includegraphics[width=0.85\textwidth]{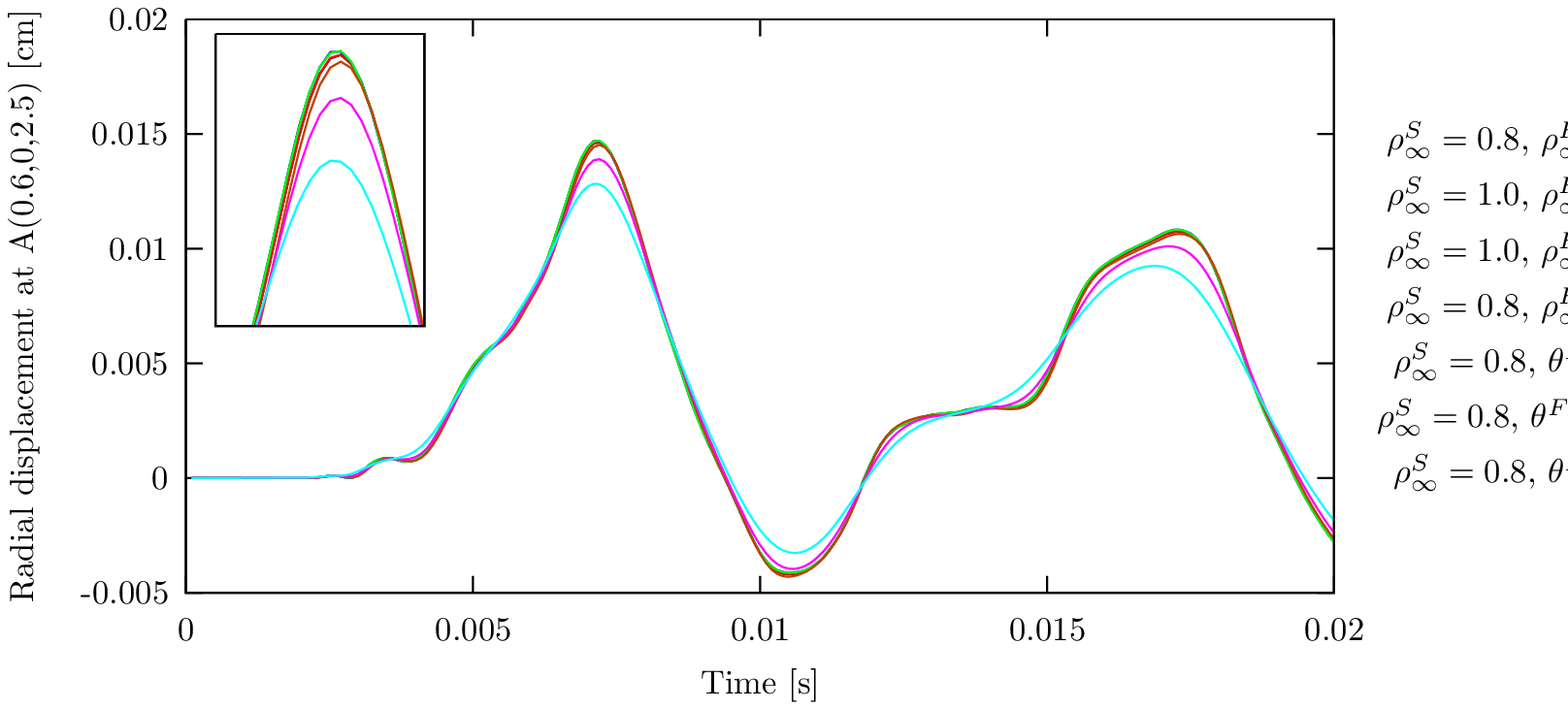}}\\
  \subfigure[fluid pressure]{\includegraphics[width=0.85\textwidth]{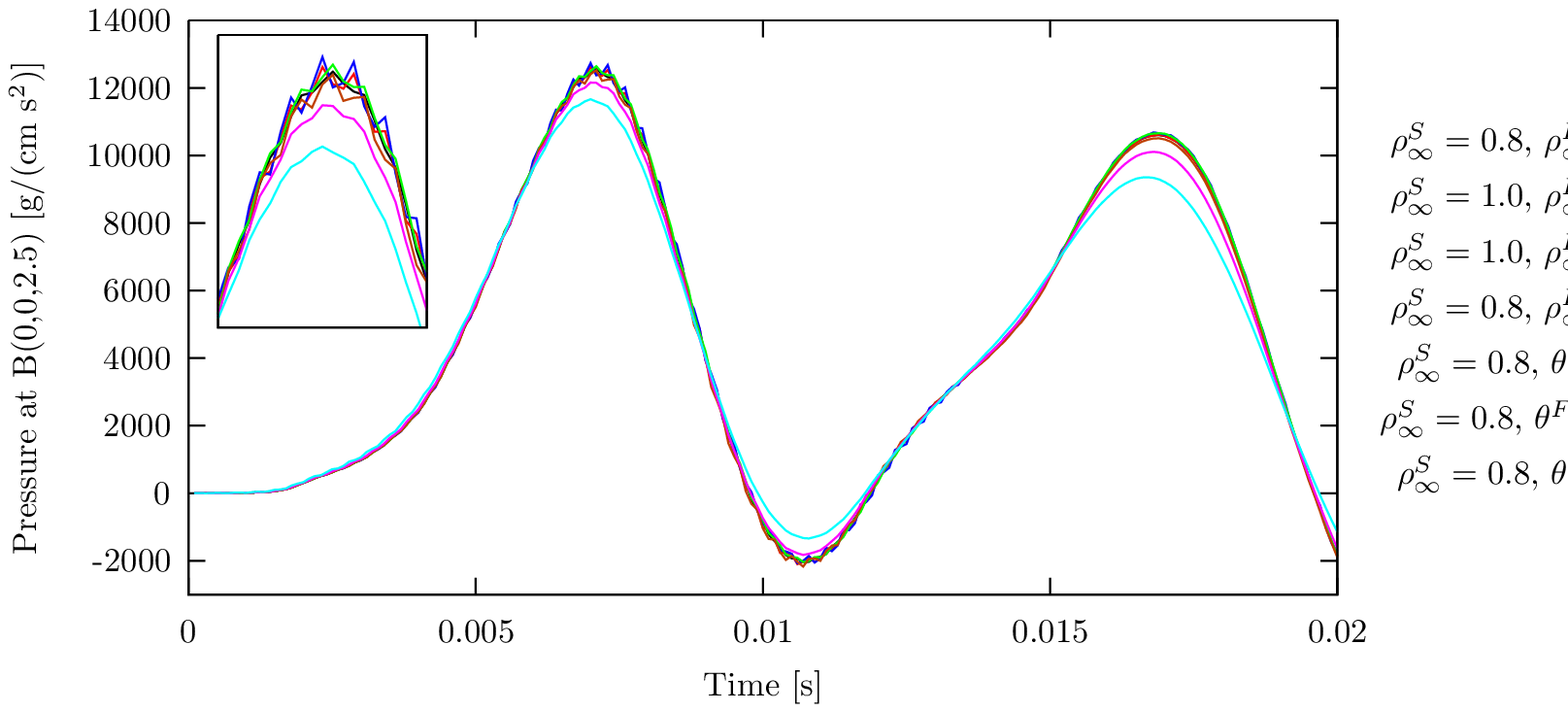}}
\caption{Various combinations of time integration schemes for pressure wave example --- For the structure field, generalized-$\alpha$ time integration with its spectral radius~$\strsub{\rho}{\infty}$ has been utilized. In the fluid field, either generalized-$\alpha$ time integration or a one-step-$\theta$ scheme have been employed, denoted by their parameters~$\flusub{\rho}{\infty}$ and~$\flu{\theta}$, respectively. Close-ups for the first peak are shown.}
\label{fig:pw_params}  
\end{center}
\vspace{-5mm}
\end{figure}

The effect of different combinations of time integration schemes with various time integration parameters in both fluid and structure field is depicted in figure~\ref{fig:pw_params}.
When both fields are discretized with generalized-$\alpha$ schemes, the influence of the actual parameter choice is rather small. The pressure field shows some fluctuations, when no or only little numerical dissipation is imposed. These fluctuations vanish the better, the more numerical dissipation is introduced into the system. However, the overall behavior of the solution is not affected by numerical dissipation. This changes dramatically, when one-step-$\theta$ time integration is utilized in the fluid field. Then, only the choice~$\flu{\theta}=0.5$ is free of numerical dissipation. The larger the value of~$\flu{\theta}$, the more numerical dissipation is involved. Again, numerical dissipation reduces the fluctuations in the pressure field. Simultaneously, the amplitudes in the displacements as well as in the pressure are reduced significantly by a larger amount of numerical dissipation, \ie the solution changes a lot. To conclude the comparison of time integration schemes, we stress that numerical dissipation in the fluid field does not only affect the fluid solution, but the solution in the structure field is also highly affected due to the interface coupling.  

\section{Concluding remarks}
\label{sec:conclusion}

A temporal consistent, mortar-based monolithic approach to large-deformation fluid-structure interaction has been proposed. It allows for both independent spatial and independent temporal discretization for both fluid and structure field. Regarding the spatial discretization, potentially non-matching grids at the fluid-structure interface are dealt with utilizing a dual mortar method. Of course, the presented framework also includes the case of conforming interface discretizations, where all mortar projection operators reduce to diagonal matrices as well as the interface constraints collapse to the trivial case of condensable point-wise constraints. For temporal discretization, both fields can be discretized using fully implicit, single-step, and single-stage time integration schemes. Thereby, the individual time integration schemes can be chosen depending on the needs of the individual fields since temporal consistency is guaranteed by the proposed method. Due to this generality, the proposed method does not impose any restrictions on the particular finite element formulations neither on fluid, ALE, nor structure field. Regarding the temporal discretization, the limitation to fully implicit, single-step, and single-stage time integration schemes does not seem harsh since they are pretty common to use. In addition, the incorporation of single-field predictors and inhomogeneous Dirichlet boundary conditions at the fluid-structure interface has been discussed.

Optimal temporal convergence rates have been shown in a simple test case with analytical solution. Furthermore, the positive effect of simple single-field predictors on the reduction of computational costs has been demonstrated. The freedom of consistently choosing different time integration schemes in fluid and structure field has been used to discuss different damping properties of the resulting algorithms.

\bibliographystyle{siam}
\bibliography{bib_mayr.bib}

\end{document}